\documentclass[12pt]{article}
\usepackage{amsfonts,amssymb,latexsym}
\newtheorem{theorem}{Theorem}[section]
\newtheorem{lemma}[theorem]{Lemma}
\newtheorem{corollary}[theorem]{Corollary}
\newtheorem{proposition}[theorem]{Proposition}
\newenvironment{proof}
{\par\addvspace{0.3cm}\noindent{\rm Proof. }}
{\nopagebreak\mbox{}\hfill $\Box$\par\addvspace{0.25cm}}
\newenvironment{proofof}[1]
{\par\addvspace{0.3cm}\noindent{\rm Proof of #1. }}
{\nopagebreak\mbox{}\hfill $\Box$\par\addvspace{0.25cm}}

\newcommand{\trace}{{\rm trace\,}}
\newcommand{\cB}{{\cal B}}
\newcommand{\cC}{{\cal C}}
\newcommand{\cL}{{\cal L}}
\newcommand{\cS}{{\cal S}}
\newcommand{\cF}{{\cal F}}
\newcommand{\cJ}{{\cal J}}
\newcommand{\cN}{{\cal N}}
\newcommand{\cW}{{\cal W}}
\newcommand{\R}{{\mathbb R}}
\newcommand{\C}{{\mathbb C}}

\newcommand{\bfH}{\mathbf{H}}
\newcommand{\bfF}{\mathbf{F}}
\renewcommand{\kappa}{\varkappa}

\newcommand{\be}{\begin{equation}}
\newcommand{\ee}{\end{equation}}

\newcommand{\bqn}{\begin{eqnarray}}
\newcommand{\eqn}{\end{eqnarray}}
\newcommand{\nn}{\nonumber}
\newcommand{\ba}{\begin{array}}
\newcommand{\ea}{\end{array}}
\newcommand{\wt}[1]{\widetilde{#1}}
\newcommand{\iv}{^{-1}}
\newcommand{\iy}{\infty}

\newcommand{\ta}{\tilde{a}}
\newcommand{\tb}{\tilde{b}}

\newcommand{\LR}[1]{L^{#1}(\R)}
\newcommand{\LRp}[1]{L^{#1}(\R_+)}
\newcommand{\Ltau}[1]{L^{#1}[0,\tau]}

\begin{document}

\date{}
\title{Asymptotics of determinants of Bessel operators}
\author{Estelle L. Basor\thanks{ebasor@calpoly.edu. 
          Supported in part by NSF Grant DMS-9970879.}\\
               Department of Mathematics\\
               California Polytechnic State University\\
               San Luis Obispo, CA 93407, USA
        \and
        Torsten Ehrhardt\thanks{tehrhard@mathematik.tu-chemnitz.de.}\\
       		Fakult\"{a}t f\"{u}r Mathematik\\
         	Technische Universit\"{a}t Chemnitz\\
         	09107 Chemnitz, Germany}
\maketitle
\begin{abstract}
For $a\in\LRp\iy\cap\LRp1$ the truncated Bessel operator $B_\tau(a)$ is
the integral operator acting on $\Ltau2$ with the kernel
$$K(x,y)=\int_0^\iy t\sqrt{xy} J_\nu(xt)J_\nu(yt)a(t)\,dt,$$
where $J_\nu$ stands for the Bessel function with $\nu>-1$.
In this paper we determine the asymptotics of the determinant
$\det(I+B_\tau(a))$ as $\tau\to\infty$ for sufficiently smooth functions
$a$ for which $a(x)\neq1$ for all $x\in[0,\iy)$.
The asymptotic formula is of the form $\det (I+B_\tau(a))\sim G^\tau E$ with certain
constants $G$ and $E$, and thus similar to the well-known
Szeg\"o-Akhiezer-Kac formula for truncated Wiener-Hopf determinants. 
\end{abstract}

%%%%%%%%%%%%%%%%%%%%%%%%%%%%%%%%%%%%%%%%%%%%%%%%%%%%%%%%%%%%%%

\section{Introduction}
\noindent

For $a\in\LRp\iy\cap\LRp1$ the {\em Bessel operator} $B(a)$ is the
integral operator acting on $\LRp2$ with the kernel
\bqn\label{f.Besselkernel}
K(x,y) &=& \int_{0}^\iy t\sqrt{xy}J_\nu(xt)J_\nu(yt)a(t)\,dt.
\eqn
Here $J_\nu$ is the Bessel function with a parameter $\nu>-1$.

For each $\tau>0$, the {\em truncated Bessel operator} $B_\tau(a)$ is the
integral operator acting on $\Ltau2$ with the same kernel
(\ref{f.Besselkernel}). Obviously, $B_\tau(a)$ can be considered as the
restriction of $B(a)$ onto $\Ltau2$, i.e., 
\bqn
B_\tau(a) &=& P_\tau B(a)P_\tau|_{\Ltau2},
\eqn
where $P_\tau$ is the projection
\bqn
P_\tau&:& f(x)\mapsto g(x)=\left\{\ba{cl}f(x)&\mbox{ for }0\le x\le\tau\\
0&\mbox{ for }x>\tau.\ea\right.
\eqn

For $a\in\LRp\iy\cap\LRp1$, the Bessel operator $B(a)$ is bounded on $\LRp2$, 
and the truncated Bessel operator $B_\tau(a)$ is a trace class operator on 
$\Ltau2$. Hence the operator determinant $\det (I+B_\tau(a))$ is well defined
for each $\tau$.
For more information about trace class operators
and related notions we refer to \cite{GK}.

In this paper we compute the asymptotics of the determinants
$\det(I+B_\tau(a))$ as $\tau\to\iy$ for certain continuous functions $a$.
Given $a\in\LRp1$, we denote by $\hat{a}$ the
cosine transform of the function $a$:
\bqn\label{f.costraf}
\hat{a}(x)&=&\frac{1}{\pi}\int_0^\iy \cos(xt)a(t)\,dt.
\eqn
A function $a$ defined on $[0,\iy)$ is said to be piecewise $C^2$ on $[0,\iy)$
if there exist $0=t_0<t_1<\dots<t_N<\iy$, $N\ge0$, such that $a$ is 
two times continuously differentiable on each of the intervals
$[0,t_1],\dots,[t_{N-1},t_N],[t_N,\iy)$, where
the derivatives at $t_0,\dots,t_N$ are considered as one-sided derivatives.

The main result of this paper is as follows.

\begin{theorem}\label{t1.1}
Let $\nu>-1$ and suppose the function $b\in\LRp\iy\cap\LRp1$ satisfies 
the following conditions:
\begin{itemize}
\item[(i)]
$b$ is continuous and piecewise $C^2$ on $[0,\iy)$, and $\lim\limits_{t\to\iy}
b(t)=0$;
\item[(ii)]
$(1+t)^{-1/2} b'(t) \in \LRp1$, $b''(t)\in \LRp1$.
\end{itemize}
Denote by $\hat{b}$ the cosine transform of $b$ and put $a=e^b-1$. Then 
\bqn\label{f.asymf}
\det(I+B_\tau(a)) &\sim& \exp\left(
\tau \hat{b}(0)-\frac{\nu}{2}b(0)+\frac{1}{2}
\int_0^\iy x(\hat{b}(x))^2\,dx\right)
\quad\mbox{ as }\tau\to\iy.
\eqn
\end{theorem}

The proof of this theorem will be given in the last section of this
paper (Section \ref{sec:6}).
Note that the assumptions of $b$ ensure that all expressions 
appearing in formula (\ref{f.asymf}) are well defined
(see also the arguments in the proof).

A result of this kind has already been established by one of the authors in
\cite{Ba} under more restrictive assumptions. There the motivation for 
considering Bessel determinants was to describe certain densities that 
occur in random matrix theory (see also \cite{BaTr}).
In comparison with \cite{Ba}, the proof of the asymptotic formula given here
will be more transparent as we also employ a new algebraic method
for the proof \cite{Eh}. In particular, we remove the quite restrictive
assumption that $||a||_{\LRp\iy}<1$, which was imposed \cite{Ba}. 
It is replaced by assumption that $1+a$ is a function that possesses a 
logarithm, which is a natural requirement 
for Szeg\"o-Akhiezer-Kac type formulas \cite{BS}.

It is notable that for particular values of $\nu$ the Bessel operator 
$B(a)$ can be written in terms of Wiener-Hopf and Hankel operators:
\bqn
B(a) &=& W(a)+H(a)\quad\mbox{ if }\nu=-1/2,\nn\\
B(a) &=& W(a)-H(a)\quad\mbox{ if }\nu=1/2.\nn
\eqn
Here we think of $a$ as a function on $\R_+$ which is extended to 
an {\em even} function on $\R$ by stipulating $a(-x)=a(x)$.
Hence in the these cases, Theorem \ref{t1.1} describes the
asymptotics of the determinants of Wiener-Hopf $+$ Hankel operators
$$ \det(I+P_\tau W(a)P_\tau\pm P_\tau H(a)P_\tau)$$
where the symbol $a$ is even. 

We should also note here that our proof requires that we establish
for general $\nu$ sufficient conditions on a function $a$ such that the Bessel 
operator $B(a)$ differs from the Wiener-Hopf operator $W(a)$ by 
a Hilbert-Schmidt operator. For $\nu = \pm 1/2$ this reduces to
the condition that the Hankel operator $H(a)$ is Hilbert-Schmidt. 
However the result for general $\nu$ is of independent interest, much more
difficult to obtain, and thus the main focus of the next section of the paper.
  
Finally, the discrete analogue of computing Toeplitz + Hankel determinants,
has been recently investigated by the authors and results have been
generalized to the case where the symbol is
discontinuous \cite{BaEh2} (see also \cite{BaEh1,BaEh3}).

%%%%%%%%%%%%%%%%%%%%%%%%%%%%%%%%%%%%%%%%%%%%%%%%%%%%%%%%%%%%%%
%%%%%%%%%%%%%%%%%%%%%%%%%%%%%%%%%%%%%%%%%%%%%%%%%%%%%%%%%%%%%%

\section{Operator theoretic preliminaries}
\label{sec:2}

In this section, we establish all general operator
theoretic facts as well particular results about Bessel operators
and Wiener-Hopf operators that we will need later on.

First of all, let us mention that Bessel operators can be defined 
for arbitrary functions $a\in\LRp\iy$. For $\nu>-1$,
let $\bfH_\nu$ denote the Hankel transform
\bqn\label{f.Htrafo}
\bfH_\nu&:&\LRp2\to\LRp2,\;
f(x)\mapsto g(x)=\int_0^\iy\sqrt{tx}J_\nu(tx)f(t)\,dt.
\eqn
It is well known that $\bfH_\nu$ is selfadjoint and unitary on $\LRp2$,
i.e., $\bfH^*_\nu=\bfH\iv_\nu=\bfH_\nu$
\cite{Ti}.
For $a\in\LRp\iy$ the Bessel operator $B(a)\in\cL(\LRp2)$ is defined by
\bqn\label{f.defbes}
B(a) &=& \bfH_\nu M(a) \bfH_\nu,
\eqn
where $M(a)$ is the multiplication operator on $\LRp2$. If $a\in\LRp\iy\cap
\LRp1$, then this definition coincides with the one given in the introduction.

{}From formula (\ref{f.defbes}) it follows immediately
that
\bqn\label{f.Bab}
B(ab) &=& B(a)B(b)
\eqn
for all $a,b\in\LRp\iy$. It is also clear that the Bessel operators are bounded
and
\bqn
\|B(a)\|_{\cL(\LRp2)} &=& \|a\|_{\LRp\iy}.
\eqn

For $a\in\LR\iy$, the two-sided Wiener-Hopf operator $W^0(a)\in\cL(\LR2)$
is defined by
\bqn\label{f.W0}
W^0(a) &=& {\bfF} M(a){\bfF}^{-1},
\eqn
where ${\bfF}:\LR2\to\LR2$ is the Fourier transform
and $M(a)$ stands here for the multiplication operator on $\LR2$.
The usual Wiener-Hopf and Hankel operators acting on $\LRp2$ are defined by
\be\label{f.defWH}
W(a) \;=\; PW^0(a)P|_{\LRp2},\qquad
H(a) \;=\; PW^0(a)JP|_{\LRp2},
\ee
where $(Pf)(x)=\chi_{\R_+}(x)f(x)$ and $(Jf)(x)=f(-x)$. For $a\in\LR\iy$,
these operators are bounded and
\be
\|W(a)\|_{\cL(\LRp2)}\;=\;\|a\|_{\LR\iy},\qquad
\|H(a)\|_{\cL(\LRp2)}\;\le\;\|a\|_{\LR\iy}.
\ee
Moreover, for $a,b\in\LR\iy$ the well known identities
\bqn\label{f.Wab}
W(ab) &=& W(a)W(b)+H(a)H(\tilde{b}),\\
H(ab) &=& W(a)H(b)+W(a)H(\tilde{b}),\label{f.Hab}
\eqn
hold, where $\tilde{b}(x)=b(-x)$. These identities are a simple consequence 
of the facts that $W^0(ab)=W^0(a)W^0(b)$, $I=P+JPJ$ and 
$JW^0(b)=W^0(\tilde{b})J$.

If $a\in\LR\iy\cap\LR1$, then $W(a)$ and $H(a)$ are integral operators on
$\LRp2$ with kernel $\hat{a}(x-y)$ and $\hat{a}(x+y)$, respectively, where
\bqn\label{f.Ftrafo}
\hat{a}(x) &=& \frac{1}{2\pi}\int_{-\iy}^\iy e^{-ixt}a(t)\,dt
\eqn
is the Fourier transform of $a$. 

We remark that if we extend $a\in\LRp1$ to an even function $a_0\in\LR1$
by stipulating $a_0(x)=a(|x|)$, $x\in\R$, then the cosine transform of $a$ coincides
with the Fourier transform of $a_0$. Therefore we will use the same notation
for the cosine transform (\ref{f.costraf}) and 
the Fourier transform (\ref{f.Ftrafo}).

In addition to the projection $P_\tau$, we define the following operators
acting on $\LRp2$:
\bqn
W_\tau&:& f(x)\mapsto g(x)=\left\{\ba{cl} f(\tau-x)&\mbox{ for } 0\le x\le\tau
\\ 0 &\mbox{ for } x>\tau,\ea\right.\\
V_\tau&:& f(x)\mapsto g(x)=\left\{\ba{cl} 0&\mbox{ for } 0\le x\le\tau
\\ f(x-\tau) &\mbox{ for } x>\tau,\ea\right.\\
V_{-\tau}&:& f(x)\mapsto g(x)=f(x+\tau),
\eqn
and $Q_\tau=I-P_\tau$. It is readily verified, that $P_\tau^2=W_\tau^2=P_\tau$,
$W_\tau P_\tau=P_\tau W_\tau=W_\tau$, $V_{-\tau}V_\tau=I$ and $V_\tau V_{-\tau}=Q_\tau$.
Moreover, the following identities hold:
\be
W_\tau W(a) W_\tau=P_\tau W(\tilde{a}) P_\tau,\quad
P_\tau W(a)V_\tau=W_\tau H(\tilde{a}),\quad
V_{-\tau}W(a)P_\tau=H(a)W_\tau.
\ee

\begin{lemma}\label{l2.0}\
\begin{itemize}
\item[(a)]
Let $a\in\LR\iy$ and $K$ be a compact operator on $\LRp2$ such that
$W(a)+K=0$. Then $a=0$ and $K=0$.
\item[(b)]
Let $a\in\LRp\iy$, $K$ be a compact operator on $\LRp2$
and $\{C_\tau\}_{\tau\in(0,\iy)}$ be a sequence of bounded operators on $\LRp2$
tending to zero in the operator norm as $\tau\to\iy$ such that
$B_{\tau}(a) + W_{\tau}KW_{\tau} + C_{\tau} = 0. $
%Is this the right statement?
Then
$a=0$, $K=0$ and $C_\tau = 0$ for all $\tau \in(0,\iy)$.
\end{itemize}
\end{lemma}
\begin{proof}
(a):\ From the equation $W(a)+K=0$ it follows that $V_{-\tau}W(a)V_\tau+V_{-\tau}KV_\tau=0$.
Since $V_{-\tau}W(a)V_\tau= W(a)$, we obtain that $W(a)=-V_{-\tau}KV_\tau$.
Observing that $V_{-\tau}\to0$ strongly as $\tau\to\iy$, we take the strong
limit of the previous equality and it follows that $W(a)=0$. Hence
$a=0$ and $K=0$.

(b):\ Since $W_\tau\to0$ weakly as $\tau\to\iy$, the operators
$W_\tau KW_\tau$ converge strongly to zero. Taking the strong limit of
$B_\tau(a)+W_\tau KW_\tau+C_\tau=0$ we obtain $B(a)=0$ because
$B_\tau(a)=P_\tau B(a)P_\tau\to B(a)$ strongly. Hence $a=0$ and
$W_\tau K W_\tau+C_\tau=0$. Multiplying with $W_\tau$ from both sides and
taking again the strong limit, we conclude that
$0=P_\tau KP_\tau+W_\tau C_\tau W_\tau\to K$. Thus $K=0$ and $C_\tau=0$.
\end{proof}

%%%%%%%%%%%%%%%%%%%%%%%%%%%%%%%%%%%%%%%%%%%%%%%%%%%%%%%

\subsection{Hilbert-Schmidt and trace class conditions}

In what follows we establish sufficient conditions for the operators
$H(a)$ and $B(a)-W(a)$ to be Hilbert-Schmidt.
Moreover, we state sufficient conditions such that $B_\tau(a)$ is a trace class
operator for each $\tau\in(0,\iy)$.

\begin{proposition}\label{p2.1}
Let $a\in \LR\iy\cap\LR1$. If
\bqn\label{f.HSnorm1}
\int_{0}^\iy x|\hat{a}(x)|^2\,dx &<&\iy,
\eqn
where $\hat{a}$ is given by (\ref{f.Ftrafo}),
then $H(a)$ is a Hilbert-Schmidt operator on $\LRp2$.
\end{proposition}
\begin{proof}
As pointed out above, $H(a)$ is an integral operator with kernel $\hat{a}(x+y)$.
This operator is Hilbert-Schmidt if and only if the following 
integral is finite, which is the square of the Hilbert-Schmidt norm of $H(a)$:
\bqn\label{f.24}
\int_0^\iy\int_0^\iy|\hat{a}(x+y)|^2\,dxdy.\nn
\eqn
This integral coincides with (\ref{f.HSnorm1}).
\end{proof}

\begin{proposition}\label{p2.2}
If $a\in\LRp\iy\cap\LRp1$, then $B_\tau(a)$ is a trace class operator
on $\Ltau2$ for each $\tau\in(0,\iy)$.
\end{proposition}
\begin{proof}
Here we make use of Mercer's Theorem \cite[Ch.III]{GK}, which reads as follows:
If $m$ is a continuous function on $[0,\tau]\times[0,\tau]$ such that
$m(s,t)=\overline{m(t,s)}$ and
\bqn
\int_0^\tau\int_0^\tau m(t,s)f(t)\overline{f(s)}\,dtds&\ge&0,
\eqn
then the (positive semi-definite) operator given by
\bqn
\Ltau2\to\Ltau2,\;f(t)\mapsto \int_0^\tau m(t,s)f(s)\,ds\nn
\eqn
is a trace class operator. This theorem shows that if $a\in\LRp\iy\cap\LRp1$
is a nonnegative function, then $B_\tau(a)$ is trace class on $\Ltau2$.
Indeed, let $m(s,t)=K(s,t)$ be the kernel (\ref{f.Besselkernel}) and
note that $K(s,t)$ is continuous since $\nu>-1$. Moreover, the above integral
(\ref{f.24}) equals
\bqn
\int_0^\iy a(t) (\bfH_\nu f)(t) \overline{(\bfH_\nu f)(t)}\,dt,\nn
\eqn
which is obviously nonnegative if $a$ is also.
Since each function in $\LRp\iy\cap\LRp1$ can be represented as the linear
combination of four nonnegative functions, the assertion follows for the
general case.
\end{proof}

The rest of this section is devoted to establishing sufficient conditions under which
the difference $B(a)-W(a)$ is a Hilbert-Schmidt operator on $\LRp2$,
where $a(-x)=a(x)$.
This result will be crucial for the considerations in the subsequent sections.

As mentioned in the introduction, it is motivated, to some extent, 
by the fact that $B(a)-W(a)=\pm H(a)$
for $\nu=\mp1/2$. In this case the desired assertion is clear from
Proposition \ref{p2.1}.

\begin{lemma}\label{l2.3}
Let $a\in\LR\iy\cap\LR1$.
Then the operators $W(a)P_1$ and $P_1W(a)$ are Hilbert-Schmidt operators
on $\LRp2$.
\end{lemma}
\begin{proof}
The operator $W(a)P_1$ is Hilbert-Schmidt if and only if 
\bqn
\int_0^1\int_0^\iy|\hat{a}(x-y)|^2\,dxdy &<&\iy,
\eqn
where $\hat{a}(x)$ is given by (\ref{f.Ftrafo}).
Obviously,
\bqn
\int_0^1\int_0^\iy|\hat{a}(x-y)|^2\,dxdy &\le&
\int_{-\iy}^\iy|\hat{a}(x)|^2\,dx
\;\;=\;\;\int_{-\iy}^\iy|a(t)|^2\,dt
\eqn
since $\hat{a}$ is the Fourier transform of $a$. This integral is finite because
$\LR\iy\cap\LR1\subset\LR2$. Hence $W(a)P_1$ is Hilbert-Schmidt.
It follows analogously that $P_1W(a)$ is a Hilbert-Schmidt operator.
\end{proof}

Let us at this point recall the following indefinite
integrals for Bessel functions:
\bqn\label{f.Bessel1}
\int tJ_\nu^2(tx)\,dt &=& \frac{t^2}{2}
\left(J_\nu^2(tx)-J_{\nu+1}(tx)J_{\nu-1}(tx)\right),\\
\label{f.Bessel2}
\int tJ_\nu(tx)J_\nu(ty)\,dt &=&
\frac{txJ_{\nu+1}(tx)J_\nu(ty)-tyJ_\nu(tx)J_{\nu+1}(ty)}{x^2-y^2}.
\eqn
They can be proved in a straightforward manner by using the recursion formulas
$J_{\nu-1}(x)+J_{\nu+1}(x)=\frac{2\nu}{x}J_\nu(x)$ and
$J_{\nu-1}(x)-J_{\nu+1}(x)=2J_\nu'(x)$.
The asymptotic behavior of $J_\nu(t)$ at zero and
at infinity is as follows:
\bqn\label{f.Bessel3}
J_\nu(t) &=& \left(\frac{t}{2}\right)^\nu\left(
\frac{1}{\Gamma(\nu+1)}+O(t^2)\right),\qquad t\to0,\\
\label{f.Bessel4}
J_\nu(t) &=& \sqrt{\frac{2}{\pi t}}\left(\cos(t-\alpha)-\sin(t-\alpha)
\frac{\nu^2-\frac{1}{4}}{t}+O(t^{-2})\right),\qquad t\to\iy,
\eqn
where $\alpha=\frac{\pi}{2}\nu+\frac{\pi}{4}$.

\begin{lemma}\label{l2.4}
Let $a\in\LRp\iy\cap\LRp1$.
Then the operators $B(a)P_1$ and $P_1B(a)$ are Hilbert-Schmidt operators
on $\LRp2$.
\end{lemma}
\begin{proof}
The operator $B(a)P_1$ is Hilbert-Schmidt if and only if 
\bqn\label{f.25}
\int_{0}^1\int_0^\iy |K(x,y)|^2\,dxdy&<&\iy,
\eqn
where $K(x,y)$ is given by (\ref{f.Besselkernel}). Another interpretation of
formula (\ref{f.Besselkernel}) is that for fixed $y$ the function $K(x,y)$
is the Hankel transform (\ref{f.Htrafo}) of the function
$a_y(t)=\sqrt{yt}J_\nu(yt)a(t)$. In other words, $K(x,y)=(\bfH_\nu a_y)(x)$.
Since $\bfH_\nu$ is an isometry on $\LRp2$, it follows that
\bqn
\int_0^\iy|K(x,y)|^2\,dx &=& \int_0^\iy|a_y(t)|^2\,dt.\nn
\eqn
Hence (\ref{f.25}) is equal to
\bqn
\int_0^1\int_0^\iy ytJ_\nu^2(yt)|a(t)|^2\,dtdy &=&
\int_0^\iy|a(t)|^2t\left(\int_0^1yJ_\nu^2(yt)\,dy\right)dt\nn\\
&=&
\frac{1}{2}\int_0^\iy|a(t)|^2t\left(J_\nu^2(t)-J_{\nu+1}(t)J_{\nu-1}(t)
\right)dt.\nn
\eqn
Here we have used (\ref{f.Bessel1}).
The term involving the Bessel functions is continuous for $t\in(0,\iy)$
and behaves at
zero like $O(t^{2\nu})$, $\nu>-1$, and at infinity like $O(t\iv)$. We can split the integral
into an integral from zero to one and an integral from one to infinity. Using the fact
that $a\in\LRp\iy\cap\LRp1$, it follows that (\ref{f.25}) is finite.
Hence $B(a)P_1$ is a Hilbert-Schmidt operator.

It can be shown analogously that $P_1B(a)$ is a Hilbert-Schmidt operator.
\end{proof}

\begin{lemma}\label{l2.5}
For each $t\in[0,\iy)$ and $x,y\in[1,\iy)$, the following integral exists:
\bqn\label{f.Kt}
K_t(x,y) &=& \int_t^\iy \left(s\sqrt{xy}J_\nu(xs)J_\nu(ys)-
\frac{2\cos(xs-\alpha)\cos(ys-\alpha)}{\pi}\right)ds.
\eqn
In particular,
\bqn\label{f.K0}
K_0(x,y) &=& -\frac{\sin(2\alpha)}{\pi(x+y)},
\eqn
and with a certain constant $C_\nu$ depending only on $\nu$ we have
\bqn\label{f.Ktnorm}
\left(\int_1^\iy\int_1^\iy |K_t(x,y)|^2\,dxdy\right)^{1/2} &\le& \frac{C_\nu}{t}
\quad\mbox{for all }t\in(0,\iy).
\eqn
\end{lemma}
\begin{proof}
{}From (\ref{f.Bessel2}), it follows that
\bqn
\lefteqn{
\int\left(s\sqrt{xy}J_\nu(xs)J_\nu(ys)-
\frac{2\cos(xs-\alpha)\cos(ys-\alpha)}{\pi}\right)ds}\hspace{4ex}\nn\\
&=& \sqrt{xy}\frac{sxJ_{\nu+1}(sx)J_\nu(sy)-syJ_\nu(sx)J_{\nu+1}(sy)}{x^2-y^2}
-\frac{\sin((x+y)s-2\alpha)}{\pi(x+y)}-\frac{\sin((x-y)s)}{\pi(x-y)}\nn\\
&=&
\sqrt{xy}\frac{sJ_{\nu+1}(sx)J_\nu(sy)+sJ_\nu(sx)J_{\nu+1}(sy)}{2(x+y)}
-\frac{\sin((x+y)s-2\alpha)}{\pi(x+y)}\nn\\
&&
\mbox{}+\sqrt{xy}\frac{sJ_{\nu+1}(sx)J_\nu(sy)-sJ_\nu(sx)J_{\nu+1}(sy)}{2(x-y)}
-\frac{\sin((x-y)s)}{\pi(x-y)}.\nn
\eqn
Using the leading term in the asymptotics (\ref{f.Bessel4}), it is easily seen
that the previous expression tends to zero as $s\to\iy$ for fixed $x,y$.
Hence the integral $K_t(x,y)$ exists for $t\in(0,\iy)$.  Using the asymptotics
(\ref{f.Bessel3}) we obtain that the terms involving the Bessel functions tend
to zero as $s\to0$ since $\nu>-1$. Hence $K_0(x,y)$ exists and equals
(\ref{f.K0}).

To be precise, the just stated assertions hold for $x\neq y$. However, if
$x=y$, we can proceed similarly by using (\ref{f.Bessel1}):
\bqn
\lefteqn{
\int\left(sxJ_\nu^2(xs)-
\frac{2\cos^2(xs-\alpha)}{\pi}\right)ds}\hspace{8ex}\nn\\
&=& \frac{xs^2}{2}\left(J_{\nu}^2(sx)-J_{\nu+1}(sx)J_{\nu-1}(sx)\right)
-\frac{\sin(2xs-2\alpha)}{2\pi x}-\frac{s}{\pi}.\nn
\eqn
Using the asymptotics (\ref{f.Bessel3}) with the first and second term
it follows that the previous expression tends to zero as $s\to\iy$.
Due to the asymptotics (\ref{f.Bessel4}) the term containing the Bessel functions
tends to zero as $s\to0$. Hence $K_t(x,x)$ exists for all $t\in[0,\iy)$ and
$K_0(x,x)$ equals (\ref{f.K0}).
 
In order to prove (\ref{f.Ktnorm}) we divide the integral into three parts:
\bqn\label{f.Ktsum}
K_t(x,y) &=& \int_t^\iy k_1(x,y;s)\,ds+
\int_t^\iy k_2(x,y;s)\,ds+\int_t^\iy k_3(x,y;s)\,ds,
\eqn
where
\bqn
k_1(x,y;s) &=& \left(\sqrt{xs}J_\nu(xs)-
\sqrt{\frac{2}{\pi}}\cos(xs-\alpha)\right)\left(\sqrt{ys}J_\nu(ys)-
\sqrt{\frac{2}{\pi}}\cos(ys-\alpha)\right),\nn\\
k_2(x,y;s) &=& \left(\sqrt{xs}J_\nu(xs)-
\sqrt{\frac{2}{\pi}}\cos(xs-\alpha)\right)\sqrt{\frac{2}{\pi}}\cos(ys-\alpha),\nn\\
k_3(x,y;s) &=& \sqrt{\frac{2}{\pi}}\cos(xs-\alpha)
\left(\sqrt{ys}J_\nu(ys)-\sqrt{\frac{2}{\pi}}\cos(ys-\alpha)\right).\nn
\eqn
Next we remark that from (\ref{f.Bessel3}) and (\ref{f.Bessel4}) it follows
that for $\nu>-1$,
\bqn\label{f.estimate1}
\left|\sqrt{z}J_\nu(z)-\sqrt{\frac{2}{\pi}}\cos(z-\alpha)\right|
&\le&\frac{\mbox{const}}{z}\qquad
\mbox{ for all }z\in(0,\iy),
\eqn
with a constant that depends only on $\nu$. Hence
$|k_1(x,y;s)|\le \mbox{const }(xys^2)\iv$, whence
\bqn
\left|\int_t^\iy k_1(x,y;s)\,ds\right|&\le&\frac{\mbox{const}}{xyt}\nn
\eqn
follows. Partial integration of the second integral in (\ref{f.Ktsum})
gives
\bqn
\int_t^\iy k_2(x,y;s)\,ds&=&
\left[\left(\sqrt{xs}J_\nu(xs)-\sqrt{\frac{2}{\pi}}\cos(xs-\alpha)\right)
\sqrt{\frac{2}{\pi}}\frac{\sin(ys-\alpha)}{y}\right]_t^\iy\nn\\
&&\hspace{-10ex}
-\int_t^\iy x\left(\frac{1}{2\sqrt{xs}}J_\nu(xs)+\sqrt{xs}J_\nu'(xs)
+\sqrt{\frac{2}{\pi}}\sin(xs-\alpha)\right)
\sqrt{\frac{2}{\pi}}\frac{\sin(ys-\alpha)}{y}\,ds.\nn
\eqn
By (\ref{f.estimate1}) the first expression is bounded by a constant times
$(xyt)\iv$. Next observe that for fixed $\nu$ and arbitrary $z\in(0,\iy)$
the following identity holds:
\bqn
\lefteqn{
\frac{1}{2\sqrt{z}}J_\nu(z)+\sqrt{z}J_{\nu}'(z) \
\;\;=\;\;
\frac{\nu+\frac{1}{2}}{\sqrt{z}}J_\nu(z)-\sqrt{z}J_{\nu+1}(z)
}\hspace{5ex}\nn\\
&=& \sqrt{\frac{2}{\pi}}\left(
\frac{\left(\nu+\frac{1}{2}\right)\cos(z-\alpha)}{z}
-\sin(z-\alpha)-\frac{\left((\nu+1)^2+\frac{1}{2}\right)\cos(z-\alpha)}{z}
+O\left(\frac{1}{z^2}\right)\right).\nn
\eqn
Like (\ref{f.estimate1}) this follows from (\ref{f.Bessel3}) and
(\ref{f.Bessel4}). We emphasize that estimate holds not just for $z\to\iy$
but also for $z\to0$, thus uniformly for all $z\in(0,\iy)$. Thus
\bqn
\int_t^\iy k_2(x,y;s)\,ds&=&
O\left(\frac{1}{xyt}\right)+\int_{t}^\iy
O\left(\frac{1}{xys^2}\right)ds+
A_\nu\int_t^\iy\frac{\cos(xs-\alpha)\sin(ys-\alpha)}{ys}\,ds\nn
\eqn
with $A_\nu=\frac{2}{\pi}(\nu^2+\nu+1)$. A similar expression can be obtained for the integral
involving $k_3(x,y;s)$. It follows that
\bqn
\lefteqn{
\int_t^\iy\left(k_2(x,y;s)+k_3(x,y;s)\right)\,ds
}\hspace{5ex}\nn\\
&=& O\left(\frac{1}{xyt}\right)
+\frac{A_\nu}{xy}\int_t^\iy\left(x\cos(xs-\alpha)\sin(ys-\alpha)+
y\sin(xs-\alpha)\cos(ys-\alpha)\right)\frac{ds}{s}.\nn
\eqn
The last integral equals
\bqn
\int_t^\iy\frac{d}{ds}\left(\sin(xs-\alpha)\sin(ys-\alpha)\right)
\frac{ds}{s}.\nn
\eqn
Another partial integration shows that this equals $O(t\iv)$.
Summarizing the previous results we can conclude that
\bqn
K_t(x,y) &=& O\left(\frac{1}{xyt}\right),\nn
\eqn
from which the desired assertion (\ref{f.Ktnorm}) follows.
\end{proof}

For $a\in\LRp\iy\cap\LRp1$, we introduce two operators acting on
$L^2[1,\iy)$. Firstly, let
\bqn
{\cal K}_a &=& Q_1(B(a)-W(a))Q_1|_{L^2[1,\iy)},
\eqn
where we stipulate $a(-x)=a(x)$, $x<0$, for the symbol of $W(a)$.
Secondly define the Hankel operator ${\cal H}_a$ as the
integral operator on $L^2[1,\iy)$ with kernel
\bqn\label{f.Ha}
{\cal H}_a(x,y)&=& -\frac{\sin(2\alpha)a(0)}{\pi(x+y)}+
\frac{1}{\pi}\int_0^\iy \cos((x+y)t-2\alpha)a(t)\,dt
\eqn
where $\alpha=\frac{\pi}{2}\nu+\frac{\pi}{4}$.

\begin{lemma}\label{l2.6}
Let $a\in\LRp\iy\cap\LRp1$ and assume that
\begin{itemize}
\item[(i)] $a$ is continuous on $[0,\iy)$;
\item[(ii)] there exist a finite number of points $0<t_1<\dots<t_N<\iy$,
$N\ge1$, such that $a$ is two times continuously differentiable on the 
interval $[0,t_1]$ and one times continuously differentiable on each of the
intervals $[t_1,t_2],\dots,[t_{N-1},t_N],[t_N,\iy)$;
\item[(iii)]
$(1+t)^{-1/2}a'(t)\in\LRp1$.
\end{itemize}
Then ${\cal K}_a-{\cal H}_a$ is a Hilbert-Schmidt operator on $L^2[1,\iy).$  
\end{lemma}
\begin{proof}
The kernel of the operator ${\cal K}_a=Q_1(B(a)-W(a))Q_1$ is given by
\bqn
{\cal K}_a(x,y) &=& \int_0^\iy \left(t\sqrt{xy}J_\nu(xt)J_\nu(yt)-
\frac{\cos(xt-yt)}{\pi}\right)a(t)\,dt.\nn
\eqn
This combined with (\ref{f.Ha}) yields that
${\cal K}_a(x,y)-{\cal H}_a(x,y)$ equals
$$
\frac{\sin(2\alpha)a(0)}{\pi(x+y)}+
\int_0^\iy\left(t\sqrt{xy}J_\nu(xt)J_\nu(yt)-
\frac{2\cos(xt-\alpha)\cos(yt-\alpha)}{\pi}\right)a(t)\,dt.
$$
In other words,
\bqn\label{f.z34}
{\cal K}_a(x,y)-{\cal H}_a(x,y)
&=&
\frac{\sin(2\alpha)a(0)}{\pi(x+y)}-\int_0^\iy
a(t)\frac{d}{dt}K_t(x,y)\,dt,
\eqn
where $K_t(x,y)$ is given by (\ref{f.Kt}).

We first consider functions $a(t)$ which are two times continuously 
differentiable on $[0,\iy)$, have compact support and satisfy $a'(0)=0$.
Notice that then $t^{-1}a'(t) \in \LRp1$.
{}From formula (\ref{f.K0}) and partial integration of (\ref{f.z34})
we obtain
\bqn\label{fz1}
{\cal K}_a(x,y)-{\cal H}_a(x,y) 
&=& \int_0^\iy a'(t)K_t(x,y)\,dt.
\eqn
Notice that $\lim\limits_{t\to\iy}K_t(x,y)=0$ and $a\in\LRp\iy$.
Equation (\ref{f.Ktnorm}) says that the integral operators on $L^2[1,\iy]$
with kernel $tK_t(x,y)$ are Hilbert-Schmidt and their Hilbert-Schmidt norm
is uniformly bounded for all $t\in(0,\iy)$. Since $t\iv a'(t)\in\LRp1$,
it follows that the operator ${\cal K}_a-{\cal H}_a$ is Hilbert-Schmidt.

Next, we assume that the function $a(t)$ satisfies the above conditions
and in addition $a(0)=0$. In this case, it is easily verified that
the following integrals
$$
\int_0^\iy\frac{|a(t)|}{t}dt,\quad
\int_{0}^{\iy}\frac{|a(t)|}{t^{3/2}(1+t^{1/2})}\,dt,\quad
\int_0^\iy\frac{|a'(t)|}{\sqrt{t}},\quad
\int_0^\iy\left|\frac{d}{dt}\left(\frac{a(t)}{t}\right)\right|dt
$$
are finite.

Since $a(0)=0$, equation (\ref{f.z34}) becomes
$$
{\cal K}_a(x,y)-{\cal H}_a(x,y) =
-\int_{0}^{\infty}a(t) \frac{d}{dt}K_{t}(x,y)\,dt .
\]
Also recall, using the notation from the proof of Lemma \ref{l2.5}, 
that we can write the function $-\frac{d}{dt}K_{t}(x,y)$ as the sum of 
three terms $k_{1}(x,y;t) + k_{2}(x,y;t)+k_{3}(x,y;t).$ We will now 
write the sum of these three operators as another sum of operators 
each of which is trace class or Hilbert-Schmidt.  To do this we recall 
that if an integral operator on $L^2[1,\iy)$ has a kernel given by
$$
K(x,y)=\int_0^\iy h_1(x,t)h_2(y,t)a(t)\,dt,
$$ 
then the trace class norm of this operator is at most
$$
\int_0^\iy|a(t)|\left(\int_1^\iy|h_1(x,t)|^2\,dx\right)^{1/2}
\left(\int_1^\iy |h_2(y,t)|^2\,dy\right)^{1/2}\,dt.
$$

Let us begin with the term $k_{1}(x,y;t)$.
{}From (\ref{f.Bessel3}) and (\ref{f.Bessel4}) along with the assumption
$\nu>-1$, we obtain
\bqn
\int_0^\iy\Big|\sqrt{z}J_\nu(z)-\sqrt{\frac{2}{\pi}}\cos(z-\alpha)\Big|^2
dz&<&\iy.\nn
\eqn 
This immediately yields that the trace class norm of the operator
given by 
$$
\int_0^\iy k_1(x,y;t)a(t) \,dt
$$ 
is bounded by a constant times $\int_0^\iy t\iv|a(t)|\,dt$.

The next term involving $k_{2}$ is a bit more complicated, 
but still follows the computation of Lemma \ref{l2.5}. 
We write 
\[
\int_{0}^{\iy}k_{2}(x,y;t)a(t)\, dt = 
\int_{0}^{\iy}\left(\sqrt{xt}J_{\nu}(xt) - \sqrt{\frac{2}{\pi}}\cos(xt 
-\alpha)\right) a(t) \sqrt{\frac{2}{\pi}}\cos(yt-\alpha)\, dt.
\] 
We use integration by parts to write this as two terms. The first term is
given by 
\[ \left[\left(\sqrt{xt}J_{\nu}(xt) - \sqrt{\frac{2}{\pi}}\cos(xt 
-\alpha)\right) a(t) \sqrt{\frac{2}{\pi}}\frac{\sin (yt -\alpha) }{y}
\right]_{t=0}^{\iy}.
\]
Since $a(t)$ is bounded at infinity, $a(t)= O(\sqrt{t})$
as $t\to0$ and $\nu>-1$, this expression is zero. 
The next term yields, from differentiating the first factor, 
two terms one of which is 
\[
\int_{0}^{\iy}\left(\sqrt{xt}J_{\nu}(xt) - \sqrt{\frac{2}{\pi}}\cos(xt 
-\alpha)\right) a'(t) \sqrt{\frac{2}{\pi}}\frac{\sin (yt -\alpha) }{y}dt.
\]
The trace norm of this operator using the same argument as above is
at most a constant times $\int_{0}^{\iy}t^{-1/2}|a'(t)|\,dt$ 
and is thus finite. So we are left with one last term
\[
\int_{0}^{\iy}
x\left(\frac{1}{2\sqrt{xt}}J_\nu(xt)+\sqrt{xt}J_\nu'(xt)
+\sqrt{\frac{2}{\pi}}\sin(xt-\alpha)\right)a(t)
\sqrt{\frac{2}{\pi}}\frac{\sin(yt-\alpha)}{y}\,dt.
\]
Following Lemma \ref{l2.5} and equation (\ref{f.Bessel4}) we rewrite 
the Bessel functions to obtain an error term of the form 
\bqn\label{f.z2}
\int_{0}^{\iy}xh(xt)a(t)\sqrt{\frac{2}{\pi}}\frac{\sin(yt -\alpha)}{y}
\,dt,
\eqn
where $h(z)$ is $O(z^{p})$ with $p=\min\{\nu -1/2,-1\}$ for small $z$
and $O(1/z^{2})$ for large $z$. The norm of the term
$$\frac{\sin(yt -\alpha)}{y}$$ in $L^2[1,\iy)$ is uniformly bounded for
$t\in[0,\iy)$ and the estimate of the norm of $xh(xt)$ in $L^2[1,\iy)$
is $O(t^{-3/2})$ as $t\to0$ and $O(t^{-2})$ as $t\to\iy$. 
Thus we have a trace norm estimate of a constant times
\[
\int_{0}^{\iy}\frac{|a(t)|}{t^{3/2}(1+t^{1/2})}\,dt,
\]
which is finite. Hence (\ref{f.z2}) has bounded trace norm.
Thus we have one remaining term
\[ A_{\nu}\int_{0}^{\iy}\frac{\cos(xt -\alpha)\sin(yt 
-\alpha)}{yt}a(t)\,dt.
\]
Now the computation for the last term $k_{3}$ is exactly the same as 
above with the variables reversed, so once again we combine these 
terms into one final integral,
\[
\frac{A_\nu}{xy}\int_0^\iy\left(x\cos(xt-\alpha)\sin(yt-\alpha)+
y\sin(xt-\alpha)\cos(yt-\alpha)\right)\frac{a(t)}{t}\,dt,
\] or 
$$
\frac{A_\nu}{xy}\int_0^\iy\frac{d}{dt}\left(\sin(xt-\alpha)\sin(yt-\alpha)
\right)\frac{a(t)}{t}\,dt. 
$$
We integrate by parts one more time so that the above is
\[
\left[\frac{A_\nu}{xy} \sin(xt-\alpha)\sin(yt-\alpha) 
\frac{a(t)}{t}\right]_{t=0}^{\iy} - \frac{A_\nu}{xy}\int_{0}^{\iy}
\sin(xt-\alpha)\sin(yt-\alpha) \left(\frac{a(t)}{t}\right)' dt.
\]
{}From our assumptions on $a(t),$ these expression are easily seen to
be $O(\frac{1}{xy})$ and are hence Hilbert-Schmidt. 

To complete the proof for arbitrary functions $a$ satisfying the assumptions
of the lemma, let $f$ be a
two times continuously differentiable function on $\R_+$ with compact support
and $f(0)=1$, $f'(0)=0$. We decompose $a=a_1+a_2$, where
$a_1(t)=a(0)f(t)$ and $a_2(t)=a(t)-a(0)f(t)$.
The function $a_1$ the fulfills the first assumed conditions, and the
function $a_2$ satisfies the second assumptions. This completes
the proof.
\end{proof}

\begin{lemma}\label{l2.7}
Let $a\in\LRp\iy\cap\LRp1$ and assume that $a$ is continuous and piecewise
$C^{2}$ on $[0,\iy)$. Assume also that
$\lim\limits_{x\to\iy}a(x)=0$ and $a''\in\LRp1$.
Then ${\cal H}_a$ is Hilbert-Schmidt.
\end{lemma}
\begin{proof}
We integrate the following integral twice by parts where $x\ge2$:
\bqn
\frac{1}{\pi}\int_0^\iy\cos(xt-2\alpha)a(t)\,dt&=&
\left[\frac{\sin(xt-2\alpha)}{\pi x}a(t)\right]_0^\iy-
\int_0^\iy\frac{\sin(xt-2\alpha)}{\pi x}a'(t)\,dt\nn\\
&=&\frac{\sin(2\alpha)a(0)}{\pi x}+\sum_{i=0}^{n-1}
\left(\left[\frac{\cos(xt-2\alpha)}{\pi x^2}a'(t)\right]_{t_i}^{t_{i+1}}
\right.\nn\\&&\left.
-\int_{t_i}^{t_{i+1}}\frac{\cos(xt-2\alpha)}{\pi x^2}a''(t)\,dt\right),\nn
\eqn
where $0=t_0<t_1<\dots t_{n-1}<t_n=\iy$ are the points where the derivatives
are discontinuous. From this it follows that
\bqn
|{\cal H}_a(x,y)|&\le& \frac{C}{(x+y)^2}.\nn
\eqn
Hence this operator is Hilbert-Schmidt.
\end{proof}

\begin{proposition}\label{p2.8}
Let $a\in\LRp\iy\cap\LRp1$ and assume that
\begin{itemize}
\item[(i)]
$a$ is continuous and piecewise $C^{2}$ on $[0,\iy)$, and 
$\lim\limits_{t\to\iy}a(t)=0$;
\item[(ii)]
$(1+t)^{-1/2}a'\in\LRp1$ and $a''\in\LRp1$.
\end{itemize}
Then the operator $B(a)-W(a)$ is Hilbert-Schmidt on $\LRp2$.
\end{proposition}
\begin{proof}
We write the operators $B(a)$ and $W(a)$ as follows:
\bqn
B(a) &=& P_1B(a)+Q_1B(a)P_1+Q_1B(a)Q_1\nn\\
W(a) &=& P_1W(a)+Q_1W(a)P_1+Q_1W(a)Q_1\nn
\eqn
Hence by Lemma \ref{l2.3} and Lemma \ref{l2.4},
\bqn
B(a)-W(a)&=&{\cal K}_a+\mbox{ Hilbert-Schmidt}.\nn
\eqn
Now it remains to apply Lemma \ref{l2.6} and Lemma \ref{l2.7}.
\end{proof}

%%%%%%%%%%%%%%%%%%%%%%%%%%%%%%%%%%%%%%%%%%%%%%%%%%%%%%%%%%%%%%
%%%%%%%%%%%%%%%%%%%%%%%%%%%%%%%%%%%%%%%%%%%%%%%%%%%%%%%%%%%%%%

\section{The algebraic approach}
\label{sec:3}

Here we follow, essentially, the method developed in \cite{Eh}. It is 
useful to point out that the next sections use Banach algebra 
techniques to compute determinants. This is done without the knowledge 
that the operator $B(a) - W(a)$ is trace class, but merely 
Hilbert-Schmidt. 

\subsection{A Banach algebra of functions}
\label{s3.1}

Let $\cS$ stand for the set of all functions $a\in\LRp\iy$ such that
the following properties are fulfilled:
\begin{itemize}
\item[(i)]
$H(a)$ is a Hilbert-Schmidt operator on $\LRp2$ where $a(-x):=a(x)$;
\item[(ii)]
$B(a)-W(a)$ is a Hilbert-Schmidt operator on $\LRp2$.
\end{itemize}
We introduce a norm in $\cS$ by stipulating
\bqn
\|a\|_{\cS} &=& \|a\|_{\LRp\iy}+\|H(a)\|_{\cC_2(\LRp2)}+
\|B(a)-W(a)\|_{\cC_2(\LRp2)}.
\eqn

\begin{proposition}
$\cS$ is a Banach algebra.
\end{proposition}
\begin{proof}
The linearity and the completeness are easy to verify. It remains to
show that $a,b\in\cS$ implies $ab\in\cS$ and that 
$\|a\|_{\cS}\|b\|_{\cS}\le \mbox{ const }\|ab\|_{\cS}$.

Indeed, let $a,b\in\cS$. Obviously $ab\in\LRp\iy$. Moreover,
\bqn
H(ab) &=& H(a)W(b)+W(a)H(b),\nn\\
B(ab)-W(ab) &=& B(a)B(b)-W(a)W(b)-H(a)H(b)\nn\\
&=& B(a)\Big(B(b)-W(b)\Big)+\Big(B(a)-W(b)\Big)W(b)-H(a)H(b).\nn
\eqn
It follows that both $H(ab)$ and $B(ab)-W(ab)$ are Hilbert
Schmidt. Moreover,
\bqn
\|H(ab)\|_{2} &\le& \|H(a)\|_{2}\|b\|_{\iy}+
\|a\|_{\iy}\|H(b)\|_{2},\nn\\
\|B(ab)-W(ab)\|_{2} &\le&
\|a\|_{\iy}\|B(b)-W(b)\|_{2}+
\|B(a)-W(b)\|_{2}\|b\|_{\iy}+
\|H(a)\|_{2}\|H(b)\|_{2}.\nn
\eqn
{}From this the norm estimate is easy to obtain.
\end{proof}

%%%%%%%%%%%%%%%%%%%%%%%%%%%%%%%%%%%%%%%%%%%%%%%%%%%%%%

\subsection{A Banach algebra of Wiener-Hopf operators}

Let $\cB$ be the set of all operators of the form
\bqn
A &=& W(a)+K
\eqn
where $a\in\cS$ and $K$ is a trace class operator on $\LRp2$.
We define in $\cB$ a norm by
\bqn
\|A\|_{\cB} &=& \|a\|_{\cS}+\|K\|_{\cC_1(\LRp2)}.
\eqn
This is definition is correct since $a$ and $K$ are uniquely determined
by the operator $A$. In fact, this is a consequence of Lemma \ref{l2.0}(a).

\begin{proposition}\label{p3.B}
$\cB$ is a Banach algebra.
\end{proposition}
\begin{proof}
The linearity and the completeness can be shown straightforwardly.
As above, we prove that $A,B\in\cB$ implies $AB\in\cB$ and a
corresponding norm estimate. In fact, let
\bqn
A\;=\;W(a)+K,\qquad B\;=\;W(b)+L,
\eqn
where $a,b\in\cS$ and $K,L\in\cC_1(\LRp2)$. Then
\bqn
AB &=& \Big(W(a)+K\Big)\Big(W(b)+L\Big)\nn\\
&=& W(ab)-H(a)H(b)+KW(b)+W(a)L+KL.\nn
\eqn
Noting that $ab\in\cS$ and $H(a)H(b)$ is trace class, 
it follows that $AB\in\cB$. Moreover,
\bqn
\|AB\|_{\cB} &\le& \|ab\|_{\cS}+\|H(a)\|_2\|H(b)\|_2+\|K\|_1\|b\|_\iy
+\|a\|_\iy\|L\|_1+\|K\|_1\|L\|_1\nn\\
&\le& \mbox{ const }\|a\|_{\cS}\|b\|_{\cS}+\|K\|_1\|b\|_{\cS}
+\|a\|_{\cS}\|L\|_1+\|K\|_1\|L\|_1.\nn
\eqn
Thus $\|AB\|_{\cB}\le\mbox{ const }\|A\|_{\cB}\|B\|_{\cB}$.
\end{proof}

%%%%%%%%%%%%%%%%%%%%%%%%%%%%%%%%%%%%%%%%%%%%%%%%%%%%%%%%%%%%%%

\subsection{A Banach algebra of sequences of Bessel operators}

We are going to introduce a Banach algebra of sequences 
$\{A_\tau\}$, $\tau\in(0,\iy)$, which contains the sequences $\{B_\tau(a)\}$ of
Bessel operators.
First of all, let $\cN$ be the set of all sequences $\{C_\tau\}$,
$\tau\in(0,\iy)$, where
\begin{itemize}
\item
$C_\tau$ is a trace class operator on $\LRp2$ for all $\tau\in(0,\iy)$;
\item
$\sup\limits_{\tau\in(0,\iy)}\|C_\tau\|_{\cC_1(\LRp2)} <\iy$;
\item
$\lim\limits_{\tau\to\iy}\|C_\tau\|_{\cC_1(\LRp2)}=0$.
\end{itemize}
Now let $\cF$ stand for the set of all sequences $\{A_\tau\}$, $\tau\in(0,\iy)$,
which are of the form
\bqn
A_\tau &=& B_\tau(a)+W_\tau KW_\tau+ C_\tau
\eqn
where $a\in\cS$, $K$ is a trace class operator on $\LRp2$ and $\{C_\tau\}\in\cN$.
For such sequences we introduce a norm by
\bqn
\left\|\{A_\tau\}\right\|_{\cF} &=&
\|a\|_{\cS}+\|K\|_{\cC_1(\LRp2)}+\sup\limits_{\tau\in(0,\iy)}
\left\|C_\tau\right\|_{\cC_1(\LRp2)}.
\eqn
This definition is correct because Lemma \ref{l2.0}(b) implies that 
for a given sequence $\{A_\tau\}$ the function $a$ and the operators $K$ and
$C_\tau$ are determined uniquely.

Moreover, for sequences $\{A_\tau^{(1)}\},\{A_\tau^{(2)}\}\in\cF$
and $\lambda^{(1)},\lambda^{(2)}\in\C$
we define algebraic operations by
\bqn
\lambda^{(1)}\{A_\tau^{(1)}\}+\lambda^{(2)}\{A_\tau^{(2)}\}=
\{\lambda^{(1)}A_\tau^{(1)}+\lambda^{(2)}A_\tau^{(2)}\},
\qquad
\{A_\tau^{(1)}\}\{A_\tau^{(2)}\}=\{A_\tau^{(1)} A_\tau^{(2)}\}.
\eqn
We will provide $\cF$ with these algebraic operations and the above norm.

\begin{proposition}
$\cF$ is a Banach algebra, and $\cN$ is a closed two-sided ideal of $\cF$.
\end{proposition}
\begin{proof}
The only non-trivial statement to prove is that 
$\{A_\tau^{(1)}\}\in\cF$ and $\{A_\tau^{(2)}\}$ implies that
$\{A_\tau^{(1)}\}\{A_\tau^{(2)}\}\in\cF$ and the corresponding norm estimate.
Let
\bqn\label{f.Atauj}
A_\tau^{(j)}&=&B_\tau(a_j)+W_\tau K_jW_\tau+C_\tau^{(j)},\quad j=1,2,
\eqn
where $a_j\in\cS$, $K_j$ is trace class and $\{C_\tau^{(j)}\}\in\cN$.
Using for brevity the notation $R(a)=B(a)-W(a)$, consider first
\bqn
P_\tau B(a_1)Q_\tau B(a_2)P_\tau &=&
P_\tau W(a_1)Q_\tau W(a_2)P_\tau+
P_\tau W(a_1)Q_\tau R(a_2)P_\tau\nn\\
&&\mbox{}+P_\tau R(a_1)Q_\tau W(a_2)P_\tau+
P_\tau R(a_1)Q_\tau R(a_2)P_\tau\nn\\
&=&
W_\tau H(a_1)H(a_2)W_\tau+
W_\tau H(a_1)V_{-\tau} R(a_2)P_\tau\nn\\
&&\mbox{}+P_\tau R(a_1)V_\tau H(a_2)W_\tau+
P_\tau R(a_1)Q_\tau R(a_2)P_\tau.\nn
\eqn
Since $H(a_j)$ and $R(a_j)$ are Hilbert-Schmidt operators and 
$V_\tau^*=V_{-\tau}$
and $Q_\tau$ converge strongly to zero on $\LRp2$ as $\tau\to\iy$, it
follows that the last three terms belong to $\cN$. Thus, it can be seen that
$\{P_\tau B(a_1)Q_\tau B(a_2)P_\tau\}\in\cF$, and the norm can be 
estimated by $\|a_1\|_{\cS}\|a_2\|_{\cS}$.

{}From this it follows that 
\bqn
B_\tau(a_1)B_\tau(a_2) &=& P_\tau B(a_1)B(a_2)P_\tau-
P_\tau B(a_1)Q_\tau B(a_2)P_\tau\nn\\
&=&P_\tau B(a_1a_2)P_\tau-W_\tau H(a_1)H(a_2)W_\tau+D_\tau^{(1)}\nn
\eqn
where $\{D_\tau^{(1)}\}\in\cN$ with 
$\|\{D_\tau^{(1)}\}\|_{\cF}\le\|a_1\|_{\cS}\|a_2\|_{\cS}$.
In particular, that $\{B_\tau(a_1)\}\{B_\tau(a_2)\}\in\cF$.

Furthermore, observe that
\bqn
B_\tau(a_1)W_\tau K_2W_\tau &=&
P_\tau W(a_1) W_\tau K_2W_\tau+P_\tau R(a_1) W_\tau K_2W_\tau\nn\\
&=& W_\tau W(a_1) P_\tau K_2W_\tau+P_\tau R(a_1) W_\tau K_2W_\tau\nn\\
&=& W_\tau W(a_1)K_2 W_\tau-W_\tau W(a_1) Q_\tau K_2W_\tau
+P_\tau R(a_1) W_\tau K_2W_\tau\nn
\eqn
where the last two terms belong to $\cN$ since
$Q_\tau\to0$ strongly and $W_\tau\to0$ weakly.
Hence 
\bqn
B_\tau(a_1)W_\tau K_2W_\tau
&=& W_\tau W(a_1)K_2 W_\tau+D_\tau^{(2)},\nn
\eqn
where $\{D_\tau^{(2)}\}\in\cN$ with 
$\|\{D_\tau^{(2)}\}\|_{\cF}\le\|a_1\|_{\cS}\|K_2\|_1$. By a similar argument
we obtain that 
\bqn
W_\tau K_1W_\tau B_\tau(a_2) &=&
W_\tau K_2W(a_2)W_\tau+D_\tau^{(3)},\nn
\eqn
where $\{D_\tau^{(3)}\}\in\cN$ with 
$\|\{D_\tau^{(3)}\}\|_{\cF}\le\|K_1\|_{1}\|a_2\|_{\cS}$. Finally,
\bqn
W_\tau K_1W_\tau W_\tau K_2W_\tau &=& W_\tau K_1K_2 W_\tau+D_\tau^{(4)},\nn
\eqn
where $\{D_\tau^{(4)}\}=\{-W_\tau K_1Q_\tau K_2W_\tau\}\in\cN$ and
$\|\{D_\tau^{(4)}\}\|_{\cF}\le\|K_1\|_{1}\|K_2\|_{1}$.

Summarizing the above we can conclude that
\bqn\label{f.A1A2}
A_\tau^{(1)}A_\tau^{(2)} &=& B_\tau(a_1a_2)+W_\tau KW_\tau+D_\tau,
\eqn
where
\bqn\label{f.K}
K&=&W(a_1)K_2+K_1W(a_2)+K_1K_2-H(a_1)H(a_2)
\eqn
and $\{D_\tau\}\in\cN$ with $\|\{D_\tau\}\|_{\cF}\le\mbox{ const }
\|\{A_\tau^{(1)}\}\|_{\cF}\|\{A_\tau^{(2)}\}\|_{\cF}$. 
Hence $\{A_\tau^{(1)}\}\{A_\tau^{(2)}\}\in\cF$. Noting that
\bqn
\|K\|_1&\le& \|a_1\|_\iy\|K_2\|_1+\|K_1\|_1\|a_2\|_\iy+\|K_1\|_1\|K_2\|_1+
\|H(a_1)\|_2\|H(a_2)\|_2,\nn\\
&\le&\left(\|a_1\|_{\cS}+\|K_1\|_1\right)\left(\|a_2\|_{\cS}+\|K_2\|_1\right),\nn
\eqn
we obtain that
\bqn
\left\|\{A_\tau^{(1)}\}\{A_\tau^{(2)}\}\right\|_{\cF}
&=&
\|a_1a_2\|_{\cS}+\|K\|_{1}+\|\{D_\tau\}\|_{\cF}\nn\\
&\le&\mbox{const }\|\{A_\tau^{(1)}\}\|_{\cF}\|\{A_\tau^{(2)}\}\|_{\cF}.\nn
\eqn
Finally, the fact that $\cN$ is an ideal of $\cF$ follows from formulas (\ref{f.A1A2})
and (\ref{f.K}) with $a_1=0$ and $K_1=0$ or $a_2=0$ and $K_2=0$, respectively.
\end{proof}

%%%%%%%%%%%%%%%%%%%%%%%%%%%%%%%%%%%%%%%%%%%%%%%%%%%%

\subsection{Banach algebra ideals and homomorphisms}

In the following proposition we introduce a Banach algebra homomorphism 
that links the Banach algebras $\cF$ and $\cB$, which have been introduced
in the previous sections. This result also shows that the 
quotient Banach algebra $\cF/\cN$ is isomorphic $\cB$. 

\begin{proposition}
The mapping $\Phi:\cF\to\cB$ defined by
$$
\Phi:\{A_\tau\}=\{B_\tau(a)+W_\tau KW_\tau+C_\tau\}\mapsto A=W(a)+K,
$$ 
is a surjective Banach algebra homomorphism with kernel $\cN$.
\end{proposition}
\begin{proof}
The linearity and surjectivity are obvious. The continuity follows from
the definition of the norms:
\bqn
\|\{A_\tau\}\|_{\cF}=\|a\|_{\cS}+\|K\|_{1}+\|\{C_\tau\}\|_{\cF}
\quad\mbox{ and }\quad
\|A\|_{\cB}=\|a\|_{\cS}+\|K\|_{1}.
\eqn
It is also clear that the kernel equals $\cN$. The only point that needs
an explanation is the multiplicativity. Assume that sequences $\{A_\tau^{(1)}\}$
and $\{A_\tau^{(2)}\}$ are given by (\ref{f.Atauj}). Then $\{A_\tau\}=
\{A_\tau^{(1)}\}\{A_\tau^{(2)}\}$ is given by formulas (\ref{f.A1A2}) and
(\ref{f.K}). Hence
\bqn
\Phi(\{A_\tau\}) &=& W(a_1a_2)+W(a_1)K_2+K_1W(a_2)+K_1K_2-H(a_1)H(a_2).\nn
\eqn
On the other hand,
\bqn
\Phi(\{A_\tau^{(j)}\}) &=& W(a_j)+K_j.\nn
\eqn
This implies that
$\Phi(\{A_\tau\})=\Phi(\{A_\tau^{(1)}\})\Phi(\{A_\tau^{(2)}\})$.
\end{proof}

Let $\cJ\subseteq\cF$ be the set
\bqn
\cJ&=&\left\{\{W_\tau KW_\tau+C_\tau\}\;:\;K\in\cC_1(\LRp2)\mbox{ and }
\{C_\tau\}\in\cN\;\right\}.
\eqn
Formulas (\ref{f.Atauj}), (\ref{f.A1A2}) and (\ref{f.K}) and the definition
of the norm in $\cF$ show that $\cJ$ is a closed two-sided ideal of
$\cF$. One can form the quotient Banach algebra $\cF/\cJ$. By
\bqn
\pi&:&\cF\to\cF/\cJ,\;\{A_\tau\}\mapsto \{A_\tau\}+\cJ,
\eqn
we denote the canonical homomorphism.

Furthermore, we introduce the linear and continuous mapping
\bqn
\Lambda&:&\cS\to\cF,\;a\mapsto\{B_\tau(a)\}.
\eqn

\begin{proposition}
The mapping $\pi\circ\Lambda:\cS\to\cF/\cJ$ is a Banach algebra isomorphism.
\end{proposition}
\begin{proof}
It is obvious that $\pi\circ\Lambda$ is linear, continuous and surjective.
In regard the multiplicativity we may refer to formulas (\ref{f.Atauj})
and (\ref{f.A1A2}) again, which can apply in order to conclude that
\bqn
B_\tau(a_1)B_\tau(a_2)=B_\tau(a_1a_2)+K+C_\tau\nn
\eqn
with a certain $K\in\cC_1(\LRp2)$ and $\{C_\tau\}\in\cN$. Finally, the
mapping is in fact an isomorphism, since for $\{A_\tau\}=\{B_\tau(a)+K+C_\tau\}$
the following norm equality holds,
\bqn
\|\{A_\tau\}+\cJ\|_{\cF/\cJ}&=&\|a\|_{\cS},\nn
\eqn
which in turn follows immediately from the definition of the norm in $\cF$.
\end{proof}

%%%%%%%%%%%%%%%%%%%%%%%%%%%%%%%%%%%%%%%%%%%%%%%%%%%%%%%%%%%%%%
%%%%%%%%%%%%%%%%%%%%%%%%%%%%%%%%%%%%%%%%%%%%%%%%%%%%%%%%%%%%%%

\section{Asymptotic analysis}
\label{sec:4}

We start with an auxiliary result concerning the exponential of an element of
$\cF$. Recall that the exponential of an element of a Banach algebra is
defined by the (absolutely convergent) sum
\bqn
e^b &=&\sum_{n=0}^\iy \frac{b^n}{n!}.
\eqn
For any Banach algebra homomorphism $\xi$ we have $\xi(e^b)=e^{\xi(b)}$.

\begin{lemma}
Let $\{A_\tau\}\in\cF$, then $\{e^{A_\tau}\}=e^{\{A_\tau\}}\in\cF$.
\end{lemma}
\begin{proof}
Obviously, $e^{\{A_\tau\}}\in\cF$. Hence $e^{\{A_\tau\}}=\{\hat{A}_\tau\}$,
where $\{\hat{A}_\tau\}$ is a sequence of operator in $\cF$. For each fixed
$\tau_0\in(0,\iy)$, the mapping $\xi_{\tau_0}:\{A_\tau\}\in\cF\mapsto
A_{\tau_0}\in\cL(\LRp2)$ is a Banach algebra homomorphism. We apply these
homomorphisms to the equation $e^{\{A_\tau\}}=\{\hat{A}_\tau\}$
and obtain that $e^{A_{\tau_0}}=\hat{A}_{\tau_0}$. Hence
$\{\hat{A}_\tau\}=\{e^{A_\tau}\}$, and the assertion is proved.
\end{proof}

The following result can be considered as a very crucial point in our
argument.

\begin{proposition}\label{p4.2}
Let $b\in\cS$ and $\{A_\tau\}=\{B_\tau(e^b)e^{-B_\tau(b)}\}$. Then there
exist an operator $K\in\cC_1(\LRp2)$ and a sequence $\{C_\tau\}\in\cN$ such
that 
\bqn\label{f.23}
\{A_\tau\}&=& \{P_\tau+W_\tau KW_\tau+C_\tau\}.
\eqn
Moreover, the operator $K$ is determined by the identity
\bqn\label{f.I+K}
I+K &=& W(e^b)e^{- W(b)}.
\eqn
\end{proposition}
\begin{proof}
First of all, the previous lemma implies that the sequence $\{A_\tau\}$ is
contained in $\cF$. Moreover, the following identity holds:
$$
\{A_\tau\} = \{B_\tau(e^b)e^{-B_\tau(b)}\}
= \{B_\tau(e^b)\}e^{-\{B_\tau(b)\}}
= \Lambda(e^b)e^{-\Lambda(b)}.
$$
Applying the canonical homomorphism $\pi$ to this identity, we obtain
$$
\pi(\{A_\tau\}) = \pi(\Lambda(e^b)e^{-\Lambda(b)})
= \pi(\Lambda(e^b))e^{-\pi(\Lambda(b))}
= \left((\pi\circ\Lambda)(e^b)\right)e^{-(\pi\circ\Lambda)(b)}.
$$
Using the fact that $\pi\circ\Lambda$ is a homomorphism, it follows that
$$
\pi(\{A_\tau\})=(\pi\circ\Lambda)(e^be^{-b})=\{P_\tau\}+\cJ.\nn
$$
Hence (\ref{f.23}) is proved. Finally we apply the homomorphism $\Phi$ to this
identity. From the definition of this homomorphism, identity (\ref{f.I+K})
follows immediately. 
\end{proof}

An operator $A$ is said to be of determinant class if $A=I+K$ where
$K$ is a trace class operator. In this case, the operator determinant
of $A$ is well defined.

A conclusion of the previous proposition is that the operator
\bqn
W(e^b)e^{-W(b)}
\eqn
is of determinant class for each $b\in\cS$. Hence the following
operator determinant
\bqn\label{f.Ea}
E[b] &=& \det W(e^b)e^{-W(b)}
\eqn
is well defined for each $b\in\cS$.
\begin{lemma}\label{l4.3}
Let $A_\tau=P_\tau+W_\tau KW_\tau+C_\tau$, where $K\in\cC_1(\LRp2)$ and
$\{C_\tau\}\in\cN$. Then
\bqn
\lim\limits_{\tau\to\iy} \det A_\tau&=&
\det(I+K).
\eqn
\end{lemma}
\begin{proof}
We remark that
$$\det A_\tau=\det(I+W_\tau KW_\tau+C_\tau)=
\det(I+P_\tau KP_\tau+W_\tau C_\tau W_\tau).$$
Since $Q_\tau^*=Q_\tau\to 0$ strongly on $\LRp2$, it follows that this equals
$\det (I+K+\widehat{C}_\tau)$ with $\widehat{C}_\tau$ tending to zero 
in the trace class norm.
\end{proof}

In order to prepare the theorem below, let us recall Proposition \ref{p2.2},
which says that for each $b\in\LRp\iy\cap\LRp1$ the operator $B_\tau(b)$ is a
trace class operator on $\Ltau2$ for each $\tau\in(0,\iy)$.
Consequently, the operator trace of $B_\tau(b)$ is well defined.  Moreover, the
following result holds.

\begin{proposition}\label{p4.extra}
For each $b\in\LRp\iy\cap\LRp1$ the operator $B_\tau(e^b)$ 
is of determinant class for each $\tau\in(0,\iy)$.
\end{proposition}
\begin{proof}
The set $\LRp\iy\cap\LRp1$ is a Banach algebra without unit element.
Upon adjoining the constant functions on $\R_+$ to $\LRp\iy\cap\LRp1$,
one obtains a Banach algebra with unit element.
{}From the series expension of $e^b$ it follows that
$e^b-1\in\LRp\iy\cap\LRp1$ whenever $b\in\LRp\iy\cap\LRp1$.
Now Proposition \ref{p2.2} implies that $B_\tau(e^b-1)$ is a trace class
operator. Thus $B_\tau(e^b)=P_\tau+B_\tau(e^b-1)$ is of determinant class
for each $\tau\in(0,\iy)$.
\end{proof}

These considerations show that all the expressions appearing in equation
(\ref{f.limit1}) below are well defined.  This equation is the asymptotic result
in which the justifications of this and the previous section culminate.

\begin{theorem}\label{t4.4}
Let $b\in\cS\cap\LRp1$. Then 
\bqn\label{f.limit1}
\lim\limits_{\tau\to\iy} 
\frac{\det B_\tau(e^b)}{e^{\trace B_\tau(b)}}
&=& E[b]
\eqn
where $E[b]$ is given by (\ref{f.Ea}).
\end{theorem}
\begin{proof}
We first use Proposition \ref{p4.2} in connection with Lemma \ref{l4.3}
and obtain that
\bqn
\lim\limits_{\tau\to\iy} \det\left(B_\tau(e^b)e^{-B_\tau(b)}\right)
&=& \det (I+K)\nn
\eqn
where $I+K$ is given by (\ref{f.I+K}). Thus $\det(I+K)$ equals $E[b]$.
Now we observe that
\bqn
\det\left(B_\tau(e^b)e^{-B_\tau(b)}\right)
&=& \left(\det B_\tau(e^b)\right)e^{-\trace B_\tau(b)},\nn
\eqn
which is a consequence of the just stated facts that $B_\tau(b)$ is
a trace class operator and $B_\tau(e^b)$ is of determinant class.
Notice that $\cS\subseteq\LRp\iy$.
\end{proof}

In order to prove the main result stated in the introduction, it is now
clear that we have to find an explicit expression for the operator determinant 
$E[b]$ and to determine the asymptotics of $\trace B_\tau(b)$ as
$\tau\to\iy$. This will be done next.

%%%%%%%%%%%%%%%%%%%%%%%%%%%%%%%%%%%%%%%%%%%%%%%%%%%%%%%%%%%%%%
%%%%%%%%%%%%%%%%%%%%%%%%%%%%%%%%%%%%%%%%%%%%%%%%%%%%%%%%%%%%%%

\section{Further calculations}
\label{sec:5}

\subsection{Trace determination}

In what follows we evaluate the asymptotics of the trace of $B_\tau(b)$ 
as $\tau\to\iy$.  The computation goes along the lines of \cite[Sect.~3]{Ba}.

\begin{proposition}\label{p5.1}
Let $b\in\LRp\iy\cap\LRp1$, and assume that $b$ is continuous and piecewise $C^1$
on $[0,\iy)$ such that $(1+t)\iv b'\in\LRp1$.  Then
\bqn
\trace B_\tau(b) &=& \frac{\tau}{\pi}\int_0^\iy b(t)\,dt-
\frac{\nu}{2}b(0)+o(1),\qquad\tau\to\iy.
\eqn
\end{proposition}
\begin{proof}
{}From formula (\ref{f.Besselkernel}) for the kernel of $B_\tau(b)$, 
it follows that
\bqn
\trace B_\tau(b) &=&
\int_0^\tau\int_0^\iy xtJ_\nu^2(xt)b(t)\,dtdx\nn\\
&=& \int_0^\iy b(t)t\int_0^\tau xJ_\nu^2(xt)\,dxdt\nn\\
&=& \frac{\tau^2}{2}\int_0^\iy b(t)t\Big(J_\nu^2(\tau t)-J_{\nu+1}(\tau t)
J_{\nu-1}(\tau t)\Big)dt\nn\\
&=& \int_0^\iy b\Big(\frac{t}{\tau}\Big)\frac{t}{2}
\Big(J_\nu^2(t)-J_{\nu+1}(t)J_{\nu-1}(t)\Big)dt.\nn
\eqn
Here (\ref{f.Bessel1}) and (\ref{f.Bessel3}) has been used. 
{}Due to the asymptotics (\ref{f.Bessel4}) we obtain
\bqn
\frac{t}{2}\Big(J_\nu^2(t)-J_{\nu+1}(t)J_{\nu-1}(t)\Big)&=&
\frac{1}{\pi}+\frac{\sin(2(t-\alpha))}{t}+O(t^{-2}),\quad
\mbox{ as }t\to\iy.\label{f.asympt1}
\eqn
Hence the following integral exists:
\bqn
F(\xi) &=& -\int_\xi^\iy\left(\frac{t}{2}\Big(J_\nu^2(t)-
J_{\nu+1}(t)J_{\nu-1}(t)\Big)-\frac{1}{\pi}\right)dt.\nn
\eqn
{}From integration by parts it follows that
\bqn
\trace B_\tau(b)-\frac{1}{\pi}\int_0^\iy b\Big(\frac{t}{\tau}\Big)\,dt
&=&
\left[b\Big(\frac{t}{\tau}\Big)F(t)\right]_{t=0}^\iy-
\frac{1}{\tau}\int_0^\iy b'\Big(\frac{t}{\tau}\Big)F(t)\,dt.\nn
\eqn
Hence 
\bqn
\trace B_\tau(b) &=&
\frac{\tau}{\pi}\int_0^\iy b(t)\,dt-b(0)F(0)-
\int_0^\iy b'(t)F(t\tau)\,dt.\nn
\eqn
The asymptotics (\ref{f.asympt1}) implies that $F(t)=O(t\iv)$ as $t\to\iy$.
Obviously, $F$ is continuous on $[0,\iy)$. Hence there exists a constant $C>0$
such that $F(t)\le C(1+t)\iv$ for all $t\in[0,\iy)$.
We estimate the last integral of the last equation as follows:
\bqn
\left|\int_0^\iy b'(t)F(t\tau)\,dt\right|
&\le&\int_0^\iy|b'(t)|\frac{C}{1+\tau t}\,dt\nn\\
&\le&
\int_0^{\frac{1}{\sqrt{\tau}}}|b'(t)|\frac{C}{1+\tau t}\,dt+
\int_{\frac{1}{\sqrt{\tau}}}^\iy|b'(t)|\frac{C}{1+\tau t}\,dt\nn\\
&\le&
C\int_0^{\frac{1}{\sqrt{\tau}}}|b'(t)|\,dt+
C\left(\sup\limits_{t\ge\frac{1}{\sqrt{\tau}}}
\frac{1+t}{1+\tau t}\right)
\int_{\frac{1}{\sqrt{\tau}}}^\iy\frac{|b'(t)|}{1+t}\,dt.\nn
\eqn
This expression tends to zero as $\tau\to\iy$. Notice that
\bqn
-F(0)&=&
\int_0^\iy\left(\frac{t}{2}\Big(J_\nu^2(t)-
J_{\nu+1}(t)J_{\nu-1}(t)\Big)-\frac{1}{\pi}\right)dt\nn\\
&=&
\int_0^\iy\left(\frac{t}{2}\Big(J_\nu^2(t)+J_{\nu+1}^2(t)\Big)-\frac{1}{\pi}
\right)dt-\nu\int_0^\iy J_{\nu+1}(t)J_\nu(t)\,dt
\;=\;-\frac{\nu}{2}.\nn
\eqn
Here we have used $J_{\nu-1}+J_{\nu+1}=\frac{2\nu}{t}J_\nu$. The first
integral equals zero, which can be seen from a direct calculation
based on (\ref{f.Bessel1}) and the asymptotics (\ref{f.Bessel3})
and (\ref{f.Bessel4}). The last integral equals $1/2$,
see, for instance, \cite[Sect.~6.512-3]{GR}.
\end{proof}

%%%%%%%%%%%%%%%%%%%%%%%%%%%%%%%%%%%%%%%%%%%%%%%%%%%%%%%%%%

\subsection{Operator determinant calculation}

In nearly the last step of our program we calculate the operator determinant
(\ref{f.Ea}) for a certain class of sufficiently smooth functions $b$.
We remark that (only) in this section we are dealing with functions
$b$ defined on $\R$ that are not necessarily even.

Let $L^2_{1/2}(\R)$ stand for the Banach space of all Lebesgue
measurable function $f$ defined on $\R$ for which
\bqn
\|f\|_{L^2_{1/2}(\R)} &=&
\left(\int_{-\iy}^\iy (1+|x|)|f(x)|^2\,dx\right)^{1/2}\;<\;\iy.
\eqn
By $C^\iy_0(\R)$ we denote the Banach algebra of all continuous functions $f$
on $\R$ for which $f(x)\to0$ as $x\to\pm\iy$, where the
multiplication is defined pointwise and the norm is the supremum norm on $\R$.

Let $\cW_0$ stand for the set of all functions $a$ defined on $\R$
which are the inverse Fourier transform of a function
$\hat{a}\in \LR1\cap L^2_{1/2}(\R)$, i.e.,
\bqn
a(x) &=& \int_{-\iy}^\iy e^{ixt}\hat{a}(t)\,dt
\eqn
(see also (\ref{f.Ftrafo})).
The norm in $\cW_0$ is defined as 
\bqn\label{f.normW0}
\|a\|_{\cW_0} &:=&
\|\hat{a}\|_{\LR1}+\|\hat{a}\|_{L^2_{1/2}(\R)}.
\eqn
A routine computation shows that the set $\LR1\cap L^2_{1/2}(\R)$ is a Banach algebra 
with the multiplication defined as the convolution on $\R$, i.e.,
\bqn
(\hat{a}\ast \hat{b})(x)=\int_{-\iy}^\iy \hat{a}(x-y)\hat{b}(y)\,dy.
\eqn
Since $\bfF\iv(\hat{a}\ast \hat{b})=(\bfF\iv\hat{a})\cdot(\bfF\iv\hat{b})$,
the set $\cW_0$ is also a Banach algebra with the multiplication
defined pointwise. Moreover,
$\cW_0$ is continuously embedded into $C^\iy_0(\R)$.

Let $\cW_0^+$ and $\cW_0^-$, respectively, be the sets of all functions 
$a\in\cW_0$, where $a=\bfF\iv \hat{a}$, for which $\hat{a}$ 
is supported on $[0,\iy)$ and $(-\iy,0]$, respectively.
It is easily seen that $\cW_0^+$ and $\cW_0^-$ are Banach subalgebras of
$\cW_0$, for which the decomposition $\cW_0=\cW_0^+\dotplus\cW_0^-$ holds.

The Banach algebras $\cW_0$ and $\cW_0^\pm$ do not contain a unit element.
The corresponding Banach algebras with unit element are obtained by adjoining
the constant functions on $\R$, i.e, 
$\cW=\C\oplus\cW_0$ and $\cW^\pm=\C\oplus\cW_0^\pm$.

Due to the fact that $\cW\subset\LR\iy$, for each $a\in\cW$
the two-sided Wiener-Hopf operator $W^0(a)$ is well defined by (\ref{f.W0}).

\begin{proposition}\label{p5.x}
Let $a\in\cW$, where $a=\alpha+\bfF\iv\hat{a}$ with $\alpha\in\C$ and
$\hat{a}\in\LR1\cap L^2_{1/2}(\R)$. Then $W^0(a)=\alpha I+K$, where 
$K$ is the integral operator on $\LR2$ 
with the kernel $k(x,y)=\hat{a}(x-y)$.
\end{proposition}
\begin{proof}
Obviously, $W^0(a)=\alpha I+W^0(\bfF\iv\hat{a})$. Let $f\in\LR1$ and denote
by $\hat{a}\ast f$ the convolution of $\hat{a}$ and $f$. Then
$\bfF\iv(\hat{a}\ast f)=(\bfF\iv\hat{a})\cdot(\bfF\iv f)$. 
Hence $\hat{a}\ast f=W^0(\bfF\iv\hat{a})f$. By an approximation argument
the same equality holds also for all functions $f\in\LR2$.
This proves the assertion.
\end{proof}

For each $a\in\cW$ the Wiener-Hopf operator $W(a)$ and the Hankel 
operator $H(a)$ are well defined in the sense of (\ref{f.defWH}).
Moreover, the previous proposition implies that under the same
hypotheses, $W(a)$ equals $\alpha I$ plus an integral operator with the
kernel $\hat{a}(x-y)$ on $\LRp2$ and $H(a)$ is an integral operator with
the kernel $\hat{a}(x+y)$ on $\LRp2$.

The following result is taken from \cite{Eh.pin}, where a formula for the
evaluation of a certain type of operator determinant has been established.  This
formula represents a generalization of Pincus' formula.  The classical Pincus'
formula has been obtained in \cite{Pi,HH} and employed in \cite{Wi76} in the
calculation of the asymptotics of Toeplitz determinants.

\begin{proposition}\label{p5.2}
Let $H$ be a Hilbert space and let $A,B\in\cL(H)$ be such that
the commutator $AB-BA$ is a trace class operator. Then 
$e^Ae^Be^{-A-B}$ is of determinant class and
\bqn
\det e^Ae^Be^{-A-B} &=& \exp \Big(\frac{1}{2}\,\trace(AB-BA)\Big).
\eqn
\end{proposition}

Next we establish an explicit formula for the operator
determinant $\det W(e^b)e^{-W(b)}$ for functions $b\in\cW$.

\begin{theorem}\label{t5.3}
Let $b\in\cW$, where $b=\beta+\bfF\iv\hat{b}$ with $\beta\in\C$ and
$\hat{b}\in \LR1\cap L^2_{1/2}(\R)$. Then
$W(e^b)e^{-W(b)}$ is of determinant class and
\be\label{f.exact}
\det W(e^b)e^{-W(b)} 
=  \exp \Big(\frac{1}{2}\trace H(b)H(\tilde{b})\Big)
=  \exp \Big(\frac{1}{2}\,\int_0^\iy 
   x\, \hat{b}(x)\hat{b}(-x)\,dx\Big). 
\ee
\end{theorem}
\begin{proof}
{}According to the remarks made after Proposition \ref{p5.x}, for 
$a\in\cW$ the Hankel operators $H(a)$ and $H(\tilde{a})$ 
(where $\tilde{a}(t)=a(t\iv)$)
are integral operators with the kernels $\hat{a}(x+y)$ and 
$\hat{a}(-x-y)$, respectively. Since $\hat{a}\in L^2_{1/2}(\R)$,
these Hankel operators are Hilbert-Schmidt operators on $\LRp2$ and the
estimate 
\bqn\label{f.63}
\|H(a)\|_{\cC_2(\LRp2)}
&\le& \|a\|_{\cW}
\eqn
\bqn\label{f.63a}
\|H(\tilde{a})\|_{\cC_2(\LRp2)}
&\le& \|a\|_{\cW}
\eqn
holds (see also the argument given in the proof of Proposition \ref{p2.1}).

Moreover, from formula (\ref{f.Wab}) it follows that
$W(ac)=W(a)W(c)$ whenever $a\in\cW_-$ or $c\in\cW_+$.
In fact, $H(a)=H(\tilde{c})=0$.
{}From this it is not difficult to conclude that
\be
W(e^{b_+})=e^{W(b_+)}\quad\mbox{ and }\quad
W(e^{b_-})=e^{W(b_-)}\label{f.67}
\ee
for all $b_+\in\cW_+$ and all $b_-\in\cW_-$.

Each function $b\in\cW$ can be decomposed into $b=b_++b_-$ with 
$b_\pm\in\cW^\pm$. We define $A=W(b_-)$ and $B=W(b_+)$. 
Using formulas (\ref{f.67}) and (\ref{f.Wab}) it follows that 
\bqn
e^Ae^B&=& e^{W(b_-)}e^{W(b_+)}=
W(e^{b_-})W(e^{b_+})=
W(e^{b_-}e^{b_+})=W(e^b).\nn
\eqn
Obviously, $e^{-A-B}=e^{-W(b_-)-W(b_+)}=e^{-W(b)}$. 
{}Again from formula (\ref{f.Wab}) it follows that
$AB-BA=W(b_-)W(b_+)-W(b_+)W(b_-)=W(b_-b_+)-W(b_+)W(b_-)=
H(b_+)H(\tilde{b}_-)=H(b)H(\tilde{b})$. 
Since both these Hankel operators
are Hilbert-Schmidt (see (\ref{f.63})),
$AB-BA$ is a trace class operator. Thus we can apply Proposition \ref{p5.2}
and obtain the first part of the desired formula.

The second part of the formula follows from the fact that the kernels of
the operators $H(b)$ and $H(\tilde{b})$ are 
given by $\hat{b}(x+y)$ and $\hat{b}(-x-y)$, respectively.
\end{proof}

For our purposes it is not enough to have formula (\ref{f.exact}) proved
for functions $b\in\cW$. In what follows we establish this formula also
for a slightly different class of functions. The proof is based on an 
approximation argument.

Let $\cS_0$ stand for the set of all functions $a\in\LR\iy$ such that 
both $H(a)$ and $H(\ta)$ are Hilbert-Schmidt operators on $\LRp2$.
We introduce a norm in $\cS_0$ as follows:
\bqn
\|a\|_{\cS_0} &=& \|a\|_{\LR\iy}+\|H(a)\|_{\cC_2(\LRp2)}+
\|H(\ta)\|_{\cC_2(\LRp2)}.
\eqn

\begin{proposition}
$\cS_0$ is a Banach algebra.
\end{proposition}
\begin{proof}
We restrict our considerations to showing that $a,b\in\cS_0$ implies
that $ab\in\cS_0$ and that 
$\|ab\|_{\cS_0}\le{\rm const}\,\|a\|_{\cS_0}\|b\|_{\cS_0}$.

In fact, if $a,b\in\cS_0$, then it follows from formula (\ref{f.Hab}) that also
$ab\in\cS_0$. Moreover,
\bqn
\|H(ab)\|_2 &=& \|a\|_\iy\|H(b)\|_2+\|H(a)\|_2\|b\|_\iy,\nn\\
\|H(\wt{ab})\|_2&=& \|a\|_\iy\|H(\tb)\|_2+\|H(\ta)\|_2\|b\|_\iy.\nn
\eqn
This implies the desired norm estimate.
\end{proof}

Let $\cB_0$ stand for the set of all operators of the form
\bqn
A&=& W(a)+K
\eqn
where $a\in\cS_0$ and $K$ is a trace class operator on $\LRp2$. We define in 
$\cB_0$ a norm by 
\bqn
\|A\|_{\cB_0} &=& \|a\|_{\cS_0}+\|K\|_{\cC_1(\LRp2)}.
\eqn
By Lemma \ref{l2.0}(a), this definition is correct.

In the same way as in the proof of Proposition \ref{p3.B} one can show
that $\cB_0$ is a Banach algebra. 
This Banach algebra possesses the ideal $\cC_1(\LRp2)$, and the canonical
homomorphism is denoted by $\pi_0:\cB_0\to\cB_0/\cC_1(\LRp2)$.
The mapping $\Lambda_0:a\in\cS_0\mapsto W(a)\in\cB_0$ is a 
continuous linear mapping. Moreover, the mapping
$\pi_0\circ\Lambda_0:\cS_0\to \cB_0/\cC_1(\LRp2)$ is a 
Banach algebra isomorphism.

\begin{proposition}\label{p.5conv}
Let $b\in\cS_0$. Then $H(b)H(\tilde{b})$ is a trace class operator and
$W(e^b)e^{-W(b)}$ is of determinant class. Moreover,
the mapping
\bqn\label{f.72}
b\in\cS_0\mapsto W(e^{b})e^{-W(b)}\in\cB_0
\eqn
is continuous.
\end{proposition}
\begin{proof}
The first assertion is obvious since both Hankel operators are Hilbert-Schmidt.

Now let $A=W(e^b)e^{-W(b)}$. Observe that $A=\Lambda_0(e^b)e^{-\Lambda_0(b)}$,
and thus $A$ belongs to $\cB_0$. Taking into account that both $\pi_0$ and 
$\pi_0\circ\Lambda_0$ are homomorphisms it follows that
$$
\pi_0(A)=\left((\pi_0\circ\Lambda_0)(e^b)\right)
e^{-(\pi_0\circ\Lambda_0)(b)}=
(\pi_0\circ\Lambda_0)(e^be^{-b})=I+\cC_1(\LRp2).
$$ 
Thus $A-I$ is a trace class operator.

The last assertion follows essentially from the fact that $\Lambda_0$ 
is continuous and that the exponential function (in $\cS_0$ and $\cB_0$)
is continuous.
\end{proof}

\begin{theorem}\label{t5.6}
Let $b\in\cS_{0}$ and assume in addition that $b\in C^\iy_0(\R)\cap\LR1$.
Then
\be\label{f.exact2}
\det W(e^b)e^{-W(b)} 
=  \exp \Big(\frac{1}{2}\trace H(b)H(\tilde{b})\Big).
\ee
\end{theorem}
\begin{proof}
Let $\hat{b}$ stand for the Fourier transform of $b$ (see (\ref{f.Ftrafo})).
Since $b\in\LR1$, the Fourier transform $\hat{b}$ is a function in
$C^\iy_0(\R)\subset\LR\iy$. 

For $\lambda\in(0,\iy)$, let $\rho_\lambda(x)=\lambda\pi\iv(1+\lambda^2x^2)\iv$
stand for the Poisson kernel, and let
$b_\lambda$ stand for the convolution of $\rho_\lambda$ with $b$, i.e.,
\bqn\label{f.poisson}
b_\lambda(x) &=& \int_{-\iy}^\iy
\frac{\lambda b(t)}{\pi(1+\lambda^2(x-t)^2)}\,dt.
\eqn
Since $\rho_\lambda$ and $b$ belong to $\LR1$, also $b_\lambda\in\LR1$.
Remark that the Fourier transform of $\rho_\lambda$ is equal to
$(\bfF\rho_\lambda)(x)=(2\pi)\iv e^{-|x|/\lambda}$. It follows that the
Fourier transform $\hat{b}_\lambda$ of $b_\lambda$ is given by
\bqn\label{f.hatbl}
\hat{b}_\lambda=\bfF(\rho_\lambda\ast b)=
2\pi(\bfF\rho_\lambda)\cdot
(\bfF b)= e^{-|x|/\lambda}\hat{b}(x).
\eqn
We see immediately that $\hat{b}_\lambda\in\LR1\cap L^2_{1/2}(\R)$. 
Thus $b_\lambda=\bfF\iv \hat{b}_\lambda$,
whence it follows that $b_\lambda\in\cW_0$.
We can apply Theorem \ref{t5.3} and obtain that
\bqn\label{f.exl}
\det W(e^{b_\lambda)})e^{-W(b_\lambda)}=
\exp\Big(\frac{1}{2}\trace H(b_\lambda) H(\tilde{b}_\lambda)\Big)
\eqn
holds for all $\lambda\in(0,\iy)$.

Another consequence of (\ref{f.hatbl}) is that 
$$ 
H(b_\lambda)=M_\lambda H(b)M_\lambda\quad\mbox{ and }\quad
H(\tb_\lambda)=M_\lambda H(\tb)M_\lambda,
$$ 
where $M_\lambda$ stands for the multiplication operator on $\LRp2$ 
with the function $e^{-x/\lambda}$, $x\in\R_+$.
Here we have to observe that the above Hankel operators are
integral operators whose kernels are given in terms of by $\hat{b}_\lambda$
and $b_\lambda$, respectively.

Since $(M_\lambda)^*=M_\lambda$ converges strongly on $\LRp2$ to the 
identity operator as $\lambda\to\iy$, it follows that
$H(b_\lambda)\to H(b)$ and $H(\tb_\lambda)\to H(\tb)$ in the Hilbert-Schmidt
norm as $\lambda\to\iy$. Hence
$$ H(b_\lambda)H(\tb_\lambda)\to H(b)H(\tb)$$
in the trace class norm as $\lambda\to\iy$.

{}From (\ref{f.poisson}) and the assumption that
$b\in C^\iy_0(\R)$ we obtain that $b_\lambda\to b$ in the norm of $\LR\iy$.
Hence, together with what has just been said,
$b_\lambda\to b$ in the norm of $\cS_0$. Now we employ the last statement of
Proposition \ref{p.5conv} and obtain that 
$$ W(e^{b_\lambda})e^{-W(b_\lambda)}\to W(e^{b})e^{-W(b)}$$
in the norm of $\cB_0$. Since both
$W(e^{b_\lambda})e^{-W(b_\lambda)}$ and $W(e^{b})e^{-W(b)}$ are operators
of the form identity plus trace class operator, it follows 
from the particular definition of the Banach algebra $\cB_0$ that
$$ W(e^{b_\lambda})e^{-W(b_\lambda)}-I\to W(e^{b})e^{-W(b)}-I$$
in the trace class norm as $\lambda\to\iy$.
By passing to the limit $\lambda\to\iy$ in (\ref{f.exl}),
the desired identity follows.
\end{proof}

The following corollary is an immediate consequence of the previous 
theorem.

\begin{corollary}\label{c5.7}
Let $b\in C^\iy_0(\R)\cap L^1(\R)$ be such that 
$\hat{b}=\bfF b\in L^2_{1/2}(\R)$. Then
\bqn
\det W(e^b)e^{-W(b)} 
&=&  \exp \Big(\frac{1}{2}\,\int_0^\iy 
     x\, \hat{b}(x)\hat{b}(-x)\,dx\Big). 
\eqn
\end{corollary}
\begin{proof}
From Proposition \ref{p2.2} it follows that the Hankel operators
$H(b)$ and $H(\tb)$ are integral operators with the kernels
$\hat{b}(x+y)$ and $\hat{b}(-x-y)$, respectively, and that both 
integral operators
are Hilbert-Schmidt. Hence $b\in\cS_0$ and $\trace H(b)H(\tb)$ 
can be expressed as the above integral.
\end{proof}
Notice that $\cW\subseteq\cS_0$. However, the example of 
$\hat{a}(x)={\rm sign}(x)e^{-|x|}$, $a(x)=2ix/(1+x^2)$ shows that
$a\in\cW_0$ but $a\notin \LR1$. Thus Corollary \ref{c5.7} does not cover
all of Theorem \ref{t5.3}.

%%%%%%%%%%%%%%%%%%%%%%%%%%%%%%%%%%%%%%%%%%%%%%%%%%%%%%%%%%%%%%%%%%

\section{Proof of the main result}
\label{sec:6}

\begin{proofof}{Theorem \ref{t1.1}}
Suppose that the function $b\in\LRp\iy\cap\LRp1$ satisfies the assumptions
(i), (ii) of Theorem \ref{t1.1}. We are going to identify $b$ with its
even continuation, $b(x)=b(-x)$.

Let $\hat{b}$ stand for the Fourier transform of the even function
$b$, i.e., the cosine transform (see (\ref{f.costraf})), and
integrate twice by parts:
\bqn
\hat{b}(x) &=& \frac{1}{\pi}\int_0^\iy\cos(xt)b(t)\,dt =
\left[\frac{\sin(xt)}{\pi x}b(t)\right]_0^\iy-
\int_0^\iy\frac{\sin(xt)}{\pi x}b'(t)\,dt\nn\\
&=&\sum_{i=0}^{n-1}
\left(\left[\frac{\cos(xt)}{\pi x^2}b'(t)\right]_{t_i}^{t_{i+1}}
-\int_{t_i}^{t_{i+1}}\frac{\cos(xt)}{\pi x^2}b''(t)\,dt\right),\nn
\eqn
where $0=t_0<t_1<\dots t_{n-1}<t_n=\iy$ are the points where the derivatives
are discontinuous. {}From this it follows that
$\hat{b}(x)=O(1/x^2)$ as $x\to\pm\iy$. Moreover, since $b\in\LRp1$, we have
$\hat{b}\in\LR\iy$. Hence we can conclude that $\hat{b}\in L^2_{1/2}(\R)$.

Firstly, this means that we are in a position to apply Corollary \ref{c5.7}.

Secondly, the Hankel operator $H(b)$ is Hilbert-Schmidt (see Proposition
\ref{p2.1}). Moreover, the assumptions (i) and (ii) also imply that
$B(b)-W(b)$ is a Hilbert-Schmidt operator (see Proposition \ref{p2.8}. Thus
from the definition of $\cS$ in Section \ref{s3.1} it follows
that $b\in\cS$. This allows us to apply Theorem \ref{t4.4}.

Finally, the assumptions (i) and (ii) also imply that the
assumptions of Proposition \ref{p5.1} are fulfilled.

Combining Theorem \ref{t4.4} with Proposition \ref{p5.1}
and Corollary \ref{c5.7} yields the desired claim of Theorem \ref{t1.1}.
\end{proofof}

%%%%%%%%%%%%%%%%%%%%%%%%%%%%%%%%%%%%%%%%%%%%%%%%%%%%%%%%%%%%%%
%%%%%%%%%%%%%%%%%%%%%%%%%%%%%%%%%%%%%%%%%%%%%%%%%%%%%%%%%%%%%%


\begin{thebibliography}{99}

\bibitem{Ba} E.L. Basor. -- Distribution functions for random variables
for ensembles of positive hermitian matrices, {\em Comm. Math. Phys.}
{\bf 188} (1997), 327--350.

\bibitem{BaEh1} E.L. Basor, T. Ehrhardt. -- On a class of 
Toeplitz + Hankel operators, {\em New York J. Math.} {\bf 5} (1999), 1--16.

\bibitem{BaEh2} E.L. Basor, T. Ehrhardt. -- Asymptotic formulas for 
determinants of a sum of finite Toeplitz and Hankel matrices,
{\em Math. Nachr.} {\bf 228} (2001), 5--45.

\bibitem{BaEh3} E.L. Basor, T. Ehrhardt. -- Asymptotic formulas for
the determinants of symmetric Toeplitz + Hankel matrices, to appear
in: {\em Operator Theory: Advances and Applications}.

\bibitem{BaTr} E.L. Basor, C.A. Tracy. --
Variance Calculations and the Bessel kernel, {\em J. Statistical Physics}
{\bf 73}, no.~2 (1993), 415--421.

\bibitem{BS} A. B\"ottcher, B. Silbermann. --
{\em Analysis of Toeplitz Operators}, Springer, Berlin, 1990.

\bibitem{Eh} T. Ehrhardt. -- A new algebraic approach to the Szeg\"o-Widom 
Limit Theorem, submitted.

\bibitem{Eh.pin} T. Ehrhardt. -- A generalization of Pincus' formula
and Toeplitz determinant identities, to appear in:
{\em Archiv der Mathematik}.

\bibitem{GK} I. Gohberg, I.M. Krein. -- {\em Introduction to the theory of
linear nonselfadjoint operators in Hilbert space}, Amer. Math. Soc.
Transl. Math. Monographs {\bf 18}, Providence, R.I., 1969.

\bibitem{GR} I.S. Gradshteyn, I.M. Ryzhik. --
Table of Integrals, Series, and Products, 5th edition, Academic Press,
San Diego, 1994.

\bibitem{HH} J.W. Helton, R.E. Howe. -- Integral operators: traces, index,
and homology, In: {\em Proceedings of the Conference in Operator Theory,
Lecture Notes in Math.}  {\bf 345}, Springer, Berlin, 1973, 141--209.

\bibitem{Me} M.L. Metha. -- {\em Random Matrices},
Academic Press, San Diego, 1991.

\bibitem{Pi} J.D. Pincus. -- On the trace of commutators in the 
algebra of operators generated by an operator with trace class
self-commutator, {\em unpublished}, 1972.

\bibitem{Ti} E.C. Titchmarsh. -- {\em Introduction to the theory of 
Fourier integrals}, Oxford, 1937.

\bibitem{Wi76} H. Widom. -- Asymptotic behavior of block Toeplitz matrices
and determinants. II, {\em Adv. Math.} {\bf 21} (1976), 1--29.

\end{thebibliography}
\end{document}